\documentclass{article}[11pt]
\usepackage[dvips]{epsfig}
\usepackage{epic}  
\usepackage{amsmath}
\usepackage{array}
\usepackage{lscape}
\usepackage{amsmath}
\usepackage{amsfonts}
\usepackage{latexsym}
\usepackage{amssymb}
\usepackage{theorem}

\usepackage{latexsym}
\usepackage{amssymb}
\usepackage{amsmath}
\usepackage{amsfonts}
\usepackage{theorem}
\usepackage{epic}  
\usepackage[usenames]{color} 

\numberwithin{equation}{section}
\theoremstyle{definition}

\newtheorem{ex}[equation]{Example}
\newtheorem{conj}[equation]{Conjecture}

\theoremstyle{plain}
\newtheorem{thm}[equation]{Theorem}
\newtheorem{prop}[equation]{Proposition}
\newtheorem{lem}[equation]{Lemma}

\newtheorem{cor}[equation]{Corollary}

\def\Young#1{\vbox{\smallskip\offinterlineskip
    \halign{&\vbox{##}\kern-\Thickness\cr #1}}}
\newdimen\Squaresize \Squaresize=16pt 
\newdimen\Thickness \Thickness=.3pt
\newdimen\Correction \Correction=7pt
\def\Gauche#1{\hbox{\vrule width \Thickness
       \vbox to \Squaresize{\vss
          \hbox to \Squaresize{\hss#1\hss}
       \vss}
    \unskip\kern\Thickness}
 \put(1,0){\line(1,0){1}}
 \put(1,1){\line(1,0){1}}
 \put(2,1){\line(0,-1){1}}   \kern-\Thickness}
\def\Vide#1{\hbox{
       \vbox to \Squaresize{\vss
          \hbox to \Squaresize{\hss#1 \hss}\vss}
    \hskip-\Correction}
   \kern-\Thickness}
\def\Droite#1{\hbox{\kern\Thickness
       \vbox to \Squaresize{\vss%
          \hbox to \Squaresize{\hss#1\hss}
       \vss}
    \unskip\vrule width \Thickness}
   \kern-\Thickness}
\def\Haut#1{\hbox{\kern-\Thickness
       \vbox to \Squaresize{\hrule height \Thickness\vss
          \hbox to \MyDim{\hss#1\hss}
       \vss}
    \unskip}
   \kern-\Thickness}
\def\Bas#1{ \hbox{\kern-\Thickness
       \vbox to \Squaresize{\vss
          \hbox to \MyDim{\hss#1\hss}

       \vss\hrule height\Thickness}
    \unskip}
   \kern-\Thickness}
\def\Carre#1{\hbox{\vrule width \Thickness
   \vbox to \Squaresize{\hrule height \Thickness\vss
      \hbox to \Squaresize{\hss#1\hss}
   \vss\hrule height\Thickness}
   \unskip\vrule width \Thickness}
   \kern-\Thickness}

\def\Box#1{\Carre{$#1$}}

\def\Sn{\cal{S}}
\def\D{\cal{D}}

\DeclareMathAlphabet{\mathdj}{U}{msb}{m}{n}
\newcommand{\R}{\ensuremath{\mathdj{R}}}
\newcommand{\N}{\ensuremath{\mathdj{N}}}
\newcommand{\Z}{\ensuremath{\mathdj{Z}}}
\newcommand{\Q}{\ensuremath{\mathdj{Q}}}
\newcommand{\Co}{\ensuremath{\mathdj{C}}}
\newcommand{\QED}{\hfill\rule{1ex}{1ex} \par\medskip}

\theoremstyle{break}

\title{\bf{Words and polynomial invariants of finite groups\\ in non-commutative variables}\\}

\author{Anouk Bergeron-Brlek, Christophe Hohlweg, Mike Zabrocki\\
{\small {\tt anouk@mathstat.yorku.ca, hohlweg.christophe@uqam.ca, zabrocki@mathstat.yorku.ca}}
}

\begin{document}
\maketitle

\begin{abstract}
 \noindent {\bf Abstract.} 
 Let $V$ be a complex vector space with basis $\{x_1,x_2,\ldots,x_n\}$ and $G$ be a finite subgroup of $GL(V)$. The tensor algebra $T(V)$ over the complex is isomorphic to the polynomials in the non-commutative variables  $x_1, x_2, \ldots, x_n$ with complex coefficients.  We want to give a combinatorial interpretation for the decomposition of $T(V)$ into simple $G$-modules. In particular, we want to study  the graded space of invariants in $T(V)$
with respect to the action of $G$. 
We give a general method for decomposing the space $T(V)$ into simple modules in terms of words in a Cayley graph of the group $G$.  To apply the method to a particular group, we require a homomorphism from a subalgebra of the group algebra into the character algebra. 
In the case of $G$ as the symmetric group, we give an example of this homomorphism from the descent algebra. When $G$ is the dihedral group, we have a realization of the character algebra as a subalgebra of the group algebra. In those two cases, we have an interpretation for the graded dimensions and the number of free generators of the algebras of invariants in terms of those words.
 \end{abstract}

\tableofcontents

\section{Introduction}
\label{sec:in}

Let $V$ be a vector space over the complex numbers $\Co$
with basis $\{ x_1, x_2, \ldots, x_n \}$ and $G$ a finite
subgroup of the general linear group $GL(V)$ of $V$, then
$$S(V) = \Co \oplus V \oplus S^2(V) \oplus S^3(V) \oplus \cdots
\simeq \Co[ x_1, x_2, \ldots, x_n]$$
is the ring of polynomials in the basis elements and
$$T(V) = \Co \oplus V \oplus V^{\otimes 2} \oplus  V^{\otimes 3} \oplus \cdots
\simeq \Co\!\left< x_1, x_2, \ldots, x_n \right>$$
is the ring of non-commutative polynomials in the basis elements
where we use the notation $S^d(V)$ to represent the $d$-fold
symmetric tensor and $V^{\otimes d} = V \otimes V \otimes \cdots \otimes V$
the $d$-fold tensor space.
We will consider the subalgebras $S(V)^G \simeq \Co[ x_1, x_2, \ldots, x_n]^G$ and
$T(V)^G\simeq \Co\left< x_1, x_2, \ldots, x_n\right>^G$
as the graded spaces of invariants with respect to the action of $G$. It is convenient to conserve the information on the dimension of each homogeneous component of degree $d$ in  the \emph{Hilbert-Poincar\'e series} $$P(S(V)^G)=\sum_{d\geq 0}{dim}~S^{d}(V)^{G}\,q^d\quad{\text and }\quad P(T(V)^G)=\sum_{d\geq 0}{dim~}(V^{\otimes d})^{G}\,q^d.$$ 

The algebra $S(V)^G$ of invariants of $G$ appears in many papers (see \cite{St} and references therein) and several algebraic tools allow us to study this algebra. MacMahon's Master theorem \cite{MM} relates the graded character
of $S(V)$ in terms of the action on $V$ by
the formula
$$\chi^{S^d(V)}(g) = \left[q^d\right]~\frac{1}{det( I - q M(g) )}$$
where $\left[q^d\right]$ represents taking the coefficient of $q^d$ in the expression
to the right and $M(g)$ is a matrix which represents the action of the group
element $g$ on a basis of $V$.
Molien's theorem \cite{Mo} allows us to calculate the graded dimensions
of the space $S(V)^G$ of invariants
$$dim~S^d(V)^G = \left[q^d\right]~
\frac{1}{|G|} \sum_{g \in G} \frac{1}{det( I - q M(g) )}~.$$
This formula alone is generally not sufficient to explain the simple
structure that we see in some of the algebras of invariants.
There is a classic result due to Chevalley \cite{Ch}
and Shephard-Todd \cite{Sh} which says that $S(V)^G$ is a free commutative
algebra if and only if $G$ is generated by pseudo-reflections (complex reflections).
Furthermore, since the groups generated by pseudo-reflections can
be classified \cite{H}, a finite list suffices to describe all the algebras of
invariants of these groups. 
A typical example for this case will be the symmetric group ${\Sn}_n$ acting on the linear span  $V={\mathcal{L}}\{ x_1, x_2, \ldots, x_n \}$ by permutation of the variables.  The invariant ring $S(V)^{{\Sn}_n}$ is isomorphic to
the ring of symmetric polynomials
in $n$ variables which is finitely generated by the $n$ invariant polynomials
$x_1^i + x_2^i + \cdots + x_n^i$ for $1 \leq i \leq n$.

In the case of the tensor algebra $T(V)$, there are analogues of  MacMahon's Master theorem and Molien's theorem \cite{DF}. The graded character can be
found in terms by what we might identify as a `master theorem'
for the tensor space,
$$\chi^{(V^{\otimes d})}(g) = tr(M(g))^d = \left[q^d\right]~\frac{1}{1- tr(M(g)) q}~.$$
The analogue of Molien's theorem \cite{DF} for the tensor
algebra says that
$$dim~(V^{\otimes d})^G = \left[q^d\right]~
\frac{1}{|G|} \sum_{g \in G} \frac{1}{1- tr(M(g)) q}~.$$
In general, we can say that $T(V)^G$ is freely
generated \cite{K,L} by an infinite set of generators (except when $G$ is scalar, {\emph i.e.} when $G$ is generated by a nonzero scalar multiple of the identity matrix) \cite{DF}.  Unlike the
case of $S(V)^G$, no simple general description of the invariants or the generators is known for large classes of groups and these
algebraic tools do not clearly show the underlying combinatorial
structure of these invariant algebras. One motivation for studying these spaces further is to find
a simple and elegant classification of the invariants of classes of groups similar to the Chevalley-Shephard-Todd result.

Our goal is to find a combinatorial interpretation for the graded dimensions and the number of free generators of these algebras of invariants which unifies their interpretations. 
The main idea is to associate to a $G$-module $V$, a subalgebra of the group algebra together with a homomorphism of algebras into the algebra of characters of $G$. 
Then we get as a consequence a combinatorial description of the graded dimensions and the free generators of the algebra $T(V)^G$ of  invariants of $G$ as paths, or words generated by a particular Cayley graph of the group $G$. More generally, we describe a general combinatorial method to decompose $T(V)$ into simple $G$-modules using those words. To compute the graded dimensions of $T(V)^G$, it then suffices to look at the multiplicity of the trivial module in $T(V)$. 

At this point, since there is not a general relation between the group algebra and the algebra of characters, we are only able to treat some examples that we decided to present in this paper and the method used gives rise to objects that are a priori not natural in that context. 
 In the case of $G$ as the symmetric group ${\Sn}_n$ on $n$ letters, the main bridge to link the words in a particular Cayley graph of ${\Sn}_n$ to the graded dimensions and the number of free generators of $T(V)^{{\Sn}_n}$ is an homomorphism from the theory of the descent algebra \cite{PR,S}. In particular, we compute the graded dimensions and the number of free generators of $T(V)^{{\Sn}_n}$ for $V$ being the geometric or the permutation module of the symmetric group ${\Sn}_n$.
When $G$ is the dihedral group or the cyclic group, we have a realization of the algebra of characters as a subalgebra of the group algebra. We then have an  interpretation for the graded dimensions and the number of free generators of the algebras of invariants in terms of words generated by a particular Cayley graph. We would be really interested to discover other examples of non-abelian groups $G$ of module $V$ for which a similar method applies. A future paper for the case of Coxeter groups of type $B$ is under construction.

When the group $G$ is generated by pseudo-reflections acting on a vector space $V$, then if $V$ is a simple $G$-module, $V$ is called the geometric module. When $G$ is the symmetric group ${\Sn}_n$ and acts on the vector space $V$ spanned by the vectors $\{ x_1, x_2,\ldots, x_n \}$ by the permutation action then $G$ is generated by pseudo-reflections, but $V$ is not a simple ${\Sn}_n$-module. The space $\Co\!\left< x_1, x_2, \ldots, x_n \right>^{{\Sn}_n}$ is known as the symmetric polynomials in non-commutative variables which was first studied by Wolf \cite{W} and more recently by Rosas-Sagan \cite{RS}. The dimension of $(V^{\otimes d})^{{\Sn}_n}$ is the number of set partitions of the numbers $\{1,2, \ldots, d\}$ into at most $n$ parts.

When $G$ is the symmetric group ${\Sn}_n$ but acting on the vector space $V$ spanned by the vectors $\{ x_1 - x_2, x_2 -x_3, \ldots, x_{n-1} -x_n \}$ (again with the permutation action on the $x_i$) then this is also a group generated by pseudo-reflections and $V$ is called the geometric module. The algebra $T(V)^{{\Sn}_n}$ of invariants is not as well understood. The graded dimensions of $T(V)^{{\Sn}_n}$ are given by the number of oscillating tableaux studied by Chauve-Goupil \cite{CG}.  This interpretation for the graded dimensions has a very different nature to that of set partitions.  By applying the results in this paper we find a combinatorial interpretation for the graded dimensions of both these spaces, and many others, which unifies the interpretations of their graded dimensions. Using the tools of the descent algebra and the Robinson-Schensted correspondence we are able to show for instance the following (see Corollary \ref{cor:set}):

Set $s_1 = (12), s_2 = (132), {\Sn}_3 = (1432), \ldots, s_{n-1} = (1\,n\,
\cdots 432)$ as elements of the symmetric group ${\Sn}_n$ in cycle notation.
The number of set partitions of the
integers $\{ 1, 2, \ldots, d \}$ into less than or equal to
$n$ parts is the number of words of length
$d$ in the alphabet $\{ e, s_1, s_2, \ldots, s_{n-1} \}$ which reduce to the identity $e$ in
the symmetric group.

The paper is organized as follows. In Section \ref{sec:ge} is described the general method used to decompose $T(V)$ into simple $G$-modules using words in a particular Cayley graph of $G$. Then  we recall in section \ref{sec:ca} the definition of a Cayley graph and present a technical lemma that we will need to link the words in a particular Cayley graph of $G$ and the decomposition of $T(V)$ into simple modules. We also present a lemma needed to link certain words and the number of free generators of the algebra $T(V)^G$ of invariants of $G$.

We consider in section \ref{sec:sy} the particular case of the symmetric group ${\Sn}_n$. Since the bridge between the words in a particular Cayley graph of ${\Sn}_n$ and the decomposition of $T(V)$ into simple modules is the theory of Solomon's descent algebra of  ${\Sn}_n$, we will first recall in section \ref{sec:so} some of its main results. In section \ref{sec:de} is proved the general theorem for the symmetric group and we then make explicit the case of $V$ being the geometric module  and the permutation one in sections \ref{sec:invess} and \ref{sec:invness} respectively.
Each of those two sections are followed by sections which contains some results and applications about the algebras $T(V)^{{\Sn}_n}$ of invariants of ${\Sn}_n$.

Finally in section \ref{sec:ot}, we apply our general method to the case of the dihedral group ${\D}_m$. We present the general theorem in section \ref{sec:genD} and in section \ref{sec:InvD} we give some results about the algebra $T(V)^{{\D}_m}$ of invariants of ${\D}_m$, when $V$ is the geometric module. 

\section{General Method}
\label{sec:ge}

Our goal is to give a combinatorial means of decomposing $T(V)$ into simple $G$-modules.
If ${\mathcal{Q}}$ is a subalgebra of the group algebra $\Co G$ and there is a surjective homomorphism from ${\mathcal{Q}}$ into the ring of characters $\Co{\text{Irr}}(G)$,  then the decomposition of $(V^{\otimes d})$ into simple modules can be computed by counting paths of length $d$ in a particular Cayley graph of $G$. 

Let us illustrate this method in details with an example. Consider the dihedral group ${\D}_3=\langle s,r\mid s^2 = r^3 = srsr=e \rangle$ with character table given in Table \ref{tabchar} and let $V^{id}$, $V^{\epsilon}$ and $V^{1}$ be the simple ${\D}_3$-modules with respectively characters
$id$, $\epsilon$ and $\chi_1$.

\begin{table}[ht!]
$$\begin{array}{|c||c|c|c|c|c|}
\hline
&\{e\}&\{r,r^2\}&\{s,rs,r^2s\}
\\\hline\hline
id&1&1&1\\\hline
\epsilon&1&-1&1\\\hline
\chi_1&2&0&-1\\\hline   
\end{array}$$
\caption{\bf Character table of ${\D}_3$.}
\label{tabchar}
 \end{table}

\begin{figure}[ht]
\begin{center}
\includegraphics[scale=0.8]{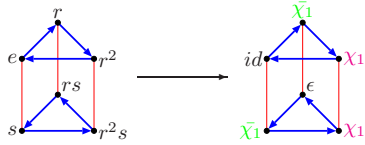}
\caption{{\bf Cayley graph of ${\D}_3$ in the generators $s$ and $r$ and decoration of its vertices with irreducible characters of ${\D}_3$.}}
\label{graphe3}
\end{center}
\end{figure}

We would like to give a combinatorial method to decompose any ${\D}_3$-module into simple modules using paths in a Cayley graph of  ${\D}_3$.
To this end, we need to find a partition of the elements of ${\D}_3$ from which we construct a subalgebra  of the group algebra $${\mathcal{Q}}={\mathcal{L}}\{e,rs,s+r, r^{2}+r^{2}s\}.$$ Then we associate each basis element of $\mathcal{Q}$ to irreducible characters of ${\D}_3$. More precisely, we define a surjective algebra morphism $\theta: {\mathcal Q}\rightarrow \Co {\rm Irr} ({\D}_3)$ by 
$$\theta(e)=id,\quad \theta(rs)=\epsilon, \quad \theta(s+r)= \bar{\chi_1}, \quad\theta( r^2+r^2s)= \chi_1,$$
where  $\bar{\chi_1}$ and $\chi_1$ represent two copies of the character $\chi_1$.  
Now if we want to decompose the module $(V^{1})^{\otimes 4}$, we consider the Cayley graph of ${\D}_3$ of Figure \ref{graphe3}, whose vertices correspond to the elements of the group and directed edges to right multiplication by the elements $s$ and $r$, since the element $r+s$ of ${\mathcal Q}$ is sent to the character $\bar{\chi_1}$ of $V^1$. Then we
 decorate the vertices with colored irreducible characters and
the multiplicity of $V^{1}$ in $(V^{1})^{\otimes 4}$, for example,  will be equal to the number of paths of length four which begin to the vertex labelled by $id$ to one labelled by $\chi_1$, plus the number of paths from the vertex labelled by $id$ to a vertex labelled by $\bar{\chi_1}$. If we pick the representatives $s$ and $r^2$ associated to the two different copies of $\chi_1$, then the multiplicities of $V^{id}$,  $V^{\epsilon}$ and  $V^{1}$ in  $({V^{1}})^{\otimes 4}$ are respectively 
$$ \begin{array}{ll}
V^{id}:&|\{ssss,rsrs,srsr\}|=3,\\ 
V^{\epsilon}:&|\{ssrs, rsss,rrsr\}|=3,\\
V^{1}:&|\{srrr,rrrs\}|+|\{ssrr,rssr,rrss\}|=5,\\
\end{array}$$
therefore $(V^{1})^{\otimes 4}$ decomposes into simple modules as
$$(V^{1})^{\otimes 4}=3\,V^{id}\oplus 3\,V^{\epsilon}\oplus 5\,V^{1}.$$
In particular, we will see later that the graded dimensions of $T(V^{1})^{{\D}_3}$ is counted by the paths from the identity to the identity in the graph of Figure \ref{graphe3}. We will also show that the number of free generators of $T(V^{1})^{{\D}_3}$ as an algebra are counted by the paths from the identity to the identity which do not cross the identity.

\section{Cayley graph of a group}
\label{sec:ca}

For the use of our purpose, let us recall the definition of a Cayley graph.
Let $G$ be a finite group and let $S\subseteq G$ be a set of group elements.
The \emph{Cayley graph} associated with $(G,S)$ is then defined as the oriented graph $\Gamma =\Gamma(G,S)$ having one vertex for each element of $G$ and the edges associated with elements in $S$. Two vertices $g_1$ and $g_2$ are joined by a directed edge associated to $s\in S$ if $g_2=g_1 s$. If the resulting Cayley graph of $G$ is connected, then the set $S$ generates $G$.

A path along the edges  of $\Gamma$ corresponds to a word in the elements in $S$. A word which \emph{reduces to $g\in G$} in $\Gamma$ is a path along the edges from the vertex corresponding to the identity to the one corresponding to the element $g$. Such a word is the reduced word corresponding to the group element $g$ with respect to the group relations. We denote by $w(g;d;\Gamma)$ the set of words of length $d$ which reduce to $g$ in $\Gamma$. We say that a word \emph{does not cross the identity} if it has no proper prefix which reduces to the identity.  

We will also consider \emph{weighted} Cayley graphs.  In other words we will associate a weight $\omega(s)$ to each element $s\in S$. We define the \emph{weight of a word} $w=s_1s_2\cdots s_r$ in the elements in $S$ to be the product of the weights of the elements in $S$, $\omega(w)=\omega(s_1)\omega(s_2)\cdots\omega(s_r)$.
To simplify the image, undirected edges represent bidirectional edges and non-labelled edges represent edges of weight one.

\begin{figure}[ht!]
\begin{center}
\includegraphics[scale=0.8]{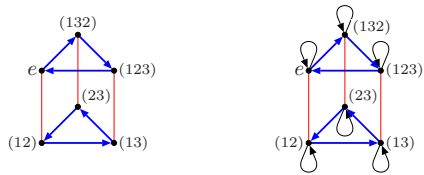}
\caption{{\bf Cayley graphs $\Gamma({\Sn}_3, \{ { {(12)}}, { {(132)}} \})$ and $\Gamma({\Sn}_3,\{e,  {(12)}, {(132)}\})$.}}
\label{graphe4}
\end{center}
\end{figure}

\begin{ex}\label{ex:cayley}
Consider the symmetric group ${\Sn}_n$ on $n$ letters with identity $e$ and permutations written in cyclic notation.
The Cayley graphs $\Gamma({\Sn}_3, \{ { {(12)}}, { {(132)}} \})$ and $\Gamma({\Sn}_3,\{e,  {(12)}, {(132)}\})$ 
are represented in Figure \ref{graphe4}.
\end{ex}

\begin{figure}[ht!]
\begin{center}
\includegraphics[scale=0.8]{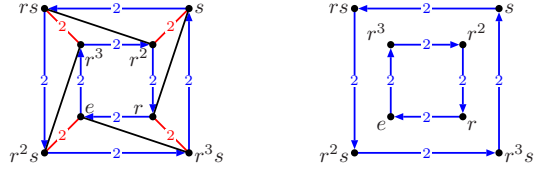}
\caption{{\bf Cayley graphs $\Gamma({\D}_4,\{{ {r^2s}}, { {r^3}},r^3s\})$ and $\Gamma({\D}_4,\{ { {r^3}}\})$ with weights $\omega(r^2s)=\omega(r^3)=2$ and $\omega(r^3s)=1$}}
\label{graphe5}
\end{center}
\end{figure}

\begin{ex}\label{ex:ca}
 Consider the dihedral group ${\D}_m=\langle { {s}}, { {r}} \mid s^2 = r^m = srsr=e \rangle$. 
The Cayley graphs $\Gamma({\D}_4,\{{ {r^2s}}, { {r^3}},r^3s\})$ and $\Gamma({\D}_4,\{ { {r^3}}\})$ with weights $\omega(r^2s)=\omega(r^3)=2$ and $\omega(r^3s)=1$ are represented in Figure \ref{graphe5}.
\end{ex}

The next key lemma will allow us to link some coefficients of an element of the group algebra to some words in a Cayley graph of $G$.

\begin{lem}\label{lemCoeff}
Let $\Gamma=\Gamma(G,\{s_1,s_2, \ldots,s_{r}\})$ be a Cayley graph of $G$ with weights $\omega(s_i)=\omega_i$. Then the coefficient of $\sigma\in G$ in the element $(\omega_1s_1+ \omega_2s_2+ \cdots +\omega_rs_{r})^d$ of the group algebra $\Co G$ equals $$\sum_{w\in w(\sigma,d; \Gamma)}\omega(w),$$ where  $w(\sigma,d; \Gamma)$ is the set of words of length $d$ which reduce to $\sigma$ in $\Gamma$.
\end{lem}
\begin{proof} Let us first prove by induction on $d$ that when the weight of each generator is one, 
the number of words of length $d$ which reduce to $\sigma$ in $\Gamma(G,\{s_1,s_2, \ldots,s_{r}\})$ is equal to the coefficient of $\sigma$ in the element $(s_1+ s_2+ \cdots +s_{r})^d$ of $\Co G$. 
Note that there is an edge from $\gamma$ to $\sigma$ in $\Gamma(G,\{s_1, s_2, \ldots,s_{r}\}) $ if and only if the coefficient of $\sigma$ in $\gamma (s_1+ s_2+ \cdots +s_{r})$ is one. 
Therefore if $d=1$, the number of paths of length one from $e$ to $\sigma$ is equal to  the coefficient of $\sigma$ in $(s_1+ s_2+ \cdots +s_{r})$.
 Let us write $s=(s_1+ s_2+ \cdots +s_{r})$ to simplify the notation.
Now if $d>1$, let $$s^d=c_{\sigma} \sigma+\sum_{\tau\in G \atop \tau\neq \sigma}c_{\tau} \tau.$$
We want to show that the number of paths of length $d$ from $e$ to $\sigma$ is equal to $c_{\sigma}$.
By induction hypothesis, the coefficient $c'_{\gamma} $ of $\gamma$ in $s^{d-1}$ is equal to the number of paths  of length $d-1$ from $e$ to $\gamma$, hence 
$$\sum_{\underbrace{\scriptstyle e\rightarrow \cdots \rightarrow \sigma}_{\rm {length}\,\,d}} 1=\sum_{\gamma\in G\atop \gamma\rightarrow \sigma}\bigg(\sum_{\underbrace{\scriptstyle e\rightarrow \cdots \rightarrow \gamma}_{\rm {length}\,\,d-1}} 1\bigg)=\sum_{\gamma\in G\atop \gamma\rightarrow \sigma} c'_{\gamma}$$
Write $s^{d-1}$ as $$s^{d-1}=\sum_{\gamma\in G\atop\gamma\rightarrow \sigma} c'_{\gamma} \gamma+\sum_{\tau\in G\atop \tau\nrightarrow \sigma}c'_{\tau} \tau.$$
Then 
$$s^{d}=s^{d-1}s=\sum_{\gamma\in G\atop \gamma\rightarrow\sigma} c'_{\gamma} \gamma s+\sum_{\tau\in G\atop\tau\nrightarrow \sigma} c'_{\tau}  \tau s.$$
Now since the coefficient of $\sigma$ in $ \gamma s$ is one if there is an edge from $\gamma$ to $\sigma$  and is zero otherwise,  
then $$\sum_{\underbrace{\scriptstyle e\rightarrow \cdots \rightarrow \sigma}_{\rm {length}\,\,d}} 1=\sum_{\gamma\in G\atop \gamma\rightarrow \sigma} c'_{\gamma}=c_{\sigma}.$$
The result then follows from the fact that a generator $s_i$ of weight $\omega_i>1$ can be seen  as a weighted-edge in the Cayley graph of $G$.
\QED\end{proof}

\begin{ex}\label{ex:motsS31}
Let us consider the Cayley graph $\Gamma=({\Sn}_3,\{(12),(132)\})$ of Figure \ref{graphe4} and set $a=(12)$ and $b=(132)$ to simplify. Then the table below shows that the coefficient of a specific element in the expansion of $(a+b)^4$  coincides with the number of words  of length three which reduce to that specific element in $\Gamma$.

\begin{center}
$$(a+b)^4={\bf{3}}\,e+{\bf{2}}\,(12)+{\bf{3}}\, (23)+{\bf{3}}\,(123)+{\bf{2}}\,(132)+{\bf{3}}\,(13)$$
\small
\begin{tabular}{|c|c|c|c|c|c|}
\hline
$e$&$(12)$&$(23)$&$(123)$&$(132)$&$(13)$\\\hline\hline
$\begin{array}[t]{c}
aaaa\\
abab\\
baba\\
\end{array}$&
$\begin{array}[t]{c}
abbb\\
bbba\\
\end{array}$&
$\begin{array}[t]{c}
aaba\\
baaa\\
bbab\\
\end{array}$&
$\begin{array}[t]{c}
aabb\\
baab\\
bbaa
\end{array}$&
$\begin{array}[t]{c}
abba\\
bbbb\\
\end{array}$&
$\begin{array}[t]{c}
aaab\\
abaa\\
babb
\end{array}$\\\hline
\end{tabular}
\normalsize
\end{center}
\end{ex}

\begin{ex}\label{ex:motsS32}
Let us consider the Cayley graph $\Gamma=({\Sn}_3,\{e,(12),(132)\})$ of Figure \ref{graphe4} and set $a=(12)$ and $b=(132)$ to simplify. Then the table below shows that the coefficient of a specific element in the expansion of $(e+a+b)^3$  coincides with the number of words  of length three which reduce to that specific element in $\Gamma$.

\begin{center}
$$(e+a+b)^3={\bf{5}}\,e+{\bf{5}}\,(12)+{\bf{4}}\, (23)+{\bf{4}}\,(123)+{\bf{5}}\,(132)+{\bf{4}}\,(13)$$
\small
\begin{tabular}{|c|c|c|c|c|c|}
\hline
$e$&$(12)$&$(23)$&$(123)$&$(132)$&$(13)$\\\hline\hline
$\begin{array}[t]{c}
bbb\\eee\\
aae\\aea\\eaa\\
\end{array}$
&
$\begin{array}[t]{c}
aaa\\bab\\
eea\\eae\\aee\\
\end{array}$
&
$\begin{array}[t]{c}
abb\\
eba\\bea\\bae\\
\end{array}$
&
$\begin{array}[t]{c}
aba\\
bbe\\beb\\ebb\\
\end{array}$
&
$\begin{array}[t]{c}
aab\\baa\\
bee\\ebe\\eeb\\
\end{array}$
&
$\begin{array}[t]{c}
bba\\
abe\\aeb\\eab\\
\end{array}$\\\hline
\end{tabular}
\normalsize
\end{center}
\end{ex}

\begin{ex}\label{ex:motsS33}
Let us consider the Cayley graph $\Gamma({\D}_4,\{{ {r^2s}}, { {r^3}},r^3s\})$ with weights $\omega(r^2s)=\omega(r^3)=2$ and $\omega(r^3s)=1$ of Figure \ref{graphe5} and set $a=r^2s$, $b=r^3$ and $c=r^3s$ to simplify. Then the table below shows that the coefficient of a specific element in the expansion of $(2a+2b+c)^2$  coincides with the sum of the weighted words of length two which reduce to that specific element in $\Gamma$.

\begin{center}
$$(2a+2b+c)^2={\bf{5}}\,e+{\bf{2}}\,r+{\bf{4}}\, r^2+{\bf{2}}\,r^3+{\bf{2}}\,s+{\bf{4}}\,rs+{\bf{2}}\,r^2s+{\bf{4}}\,r^3s$$
\small
\begin{tabular}{|c|c|c|c|c|c|c|c|}
\hline
$e$&$r$&$r^2$&$r^3$&$s$&$rs$&$r^2s$&$r^3s$\\\hline\hline
$
\omega(aa)+
\omega(cc)
$&
$\omega(ca)
$&
$
\omega(bb)
$&
$
\omega(ab)
$&
$
\omega(cb)$&
$
\omega(ba)$
&
$
\omega(bc)$
&
$
\omega(ab)$\\\hline
$5$&$2$&$4$&$2$&$2$&$4$&$2$&$4$\\\hline
\end{tabular}
\normalsize
\end{center}
\end{ex}

The next general lemma will be used to link the free generators of the invariant algebra of $G$ to words in a Cayley graph of $G$.

\begin{lem}\label{lem:gen}
Consider the Cayley graph $\Gamma$ of the group $G$ and the generating series 
$A_{\sigma}(q)$ counting the number of words which reduce to the element $\sigma\in G$  in $\Gamma$ and  $B_{\sigma}(q)$ counting the number of words which reduce to the element $\sigma\in G$ in $\Gamma$ without crossing the identity. Then we have the relation
$$A_{\sigma}(q)=\frac{1}{1-B_e(q)}B_{\sigma}(q)\quad\text{ and }\quad A_{e}(q)=\frac{1}{1-B_e(q)}.$$
\end{lem}
\begin{proof}
A path in $\Gamma$ from the identity vertex  to the vertex $\sigma\in G$  can be seen as a concatenation of paths from the identity to the identity which do not cross the identity,  and  a path from the identity to $\sigma$ which do not cross the identity
$$e\rightarrow \cdots \rightarrow e\mid e \rightarrow \cdots \rightarrow e\mid \cdots \mid e\rightarrow \cdots \rightarrow \sigma.$$
In terms of $A_{\sigma}(q)$  and $B_{\sigma}(q)$, this translates into
\begin{align*}
A_{\sigma}(q)&=B_{\sigma}(q)+B_{e}(q)B_{\sigma}(q)+(B_{e}(q))^2B_{\sigma}(q)+(B_{e}(q))^3B_{\sigma}(q)+\cdots\\
&=(1+B_{e}(q)+(B_{e}(q))^2+(B_{e}(q))^3+\cdots)B_{\sigma}(q)\\
&=\frac{1}{1-B_e(q)}B_{\sigma}(q).
\end{align*}
\QED\end{proof}

\section{Symmetric group ${\Sn}_n$}
\label{sec:sy}
We will present in this section the case of the symmetric group ${\Sn}_n$ on $n$ letters by giving  a combinatorial way to decompose the tensor algebra on $V$ into simple modules, where $V$ is any ${\Sn}_n$-module. 
We will then make explicit the cases of $V$ being the geometric module and $V$ being the permutation one, by means of words in a particular Cayley graphs of ${\Sn}_n$. We will also give a combinatorial way to compute the graded dimensions of the space $T(V)^{{\Sn}_n}$ of invariants of ${\Sn}_n$,  which is the multiplicity of the trivial module in the decomposition of $T(V)$ into simple modules, and give a description for the number of free generators of $T(V)^{{\Sn}_n}$ as an algebra. But first let us recall some definition and the theory of the descent algebra  which is the bridge between those words and the decompositon of $T(V)$. 

\subsection{Partitions and  tableaux}\label{sec:Partition}
\label{sec:pa}

\begin{figure}[ht]
\begin{center}
\includegraphics[scale=0.8]{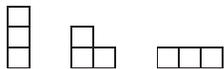}
\caption{{\bf Ferrers diagrams of $(1,1,1)$, $(2,1)$ and $(3)$.}}
\label{partages}
\end{center}
\end{figure}

\begin{figure}[ht]\begin{center}
\includegraphics[scale=0.8]{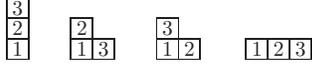}
\caption{{\bf Standard tableau with three squares.}}
\label{tableaux}
\end{center}
\end{figure}

To fix the notation, recall the definition of a partition.
A \emph{partition} $\lambda$ of a positive integer $n$ is a decreasing sequence
$\lambda_1\ge \lambda_2\ge \ldots\ge\lambda_\ell > 0$
of positive integers such that
$n=|\lambda|=\lambda_1+\lambda_2+\ldots+\lambda_\ell.$
We will write
$\lambda=(\lambda_1,\lambda_2,\ldots \lambda_\ell)\vdash n$. 
It is natural to represent a partition by a diagram. The \emph{Ferrers diagram} of a partition  $\lambda=(\lambda_1,\lambda_2,\ldots \lambda_\ell)$ is the finite subset  $\lambda=\{(a,b)\,|\,0\leq a\leq \ell-1{\text{ and }} 0\leq b\leq
\lambda_{a+1}-1\}$ of $\N\times\N$.
Visually, each element of $\lambda$ corresponds to the bottom left corner of a square of dimension $1\times 1$ in $\N\times\N$.
\begin{ex}The partitions of three are $(1,1,1)$, $(2,1)$ and $(3)$
and their Ferrers diagrams are represented in Figure \ref{partages}.
\end{ex}
A \emph{tableau} of shape $\lambda\vdash n$, denoted $shape(t)=\lambda$, with values in $T=\{1,2,\ldots,n\}$ is a function 
$
t:\lambda\longrightarrow T.
$
We can visualize it with filling each square  $c$ of a Ferrers diagram $\lambda$ with the value $t(c)$. In a natural way, a tableau is said to be injective if the function $t$ is injective. A tableau is said to be \emph{standard} if it is injective and its entries form an increasing sequence along each row and along each column. We will denote by $STab_n$ the set of standard tableau with $n$ squares.
\begin{ex} The set $STab_3$ contains the four standard tableaux of Figure \ref{tableaux}.
\end{ex}
The  \emph{Robinson-Schensted correspondence} \cite{Ro,Sc} is a bijection between the elements $\sigma$ of the symmetric group ${\Sn}_n$ and pairs $(P(\sigma),Q(\sigma))$ of standard tableaux of the same shape, where $P(\sigma)$ is the insertion tableau and $Q(\sigma)$ the recording tableau. 

\subsection{Simple ${\Sn}_n$-modules}
\label{sec:re}
Since the conjugacy classes in ${\Sn}_n$ are in bijection with the partitions of $n$, it is natural to index  the simple ${\Sn}_n$-modules by  the partitions $\lambda$ of $n$ and we will denote them by $V^{\lambda}$. In particular, the simple ${\Sn}_n$-module $V^{(n)}$ is the trivial one and $V^{(n-1,1)}$ the geometric one.
Let us consider the $\Q$-linear span $V={\mathcal{L}}\{x_1,x_2,\ldots, x_{n}\}$ on which ${\Sn}_n$ acts by permuting the coordinates. Then we have
$$V={\mathcal{L}}\{x_1+x_2+x_3+\cdots+x_{n}\}\oplus {\mathcal{L}}\{x_1-x_2,x_2-x_3,\ldots,x_{n-1}-x_{n}\}$$ so the decomposition of the permutation module $V$ into simple ${\Sn}_n$-modules is $V=V^{(n)}\oplus V^{(n-1,1)}.$ Let $X_n$ denote the set of variables $x_1,x_2,\ldots, x_n$. If we identify $T(V)$ with $\Q\langle X_{n}\rangle$, then $$T(V^{(n-1,1)})\simeq\Q\langle X_{n}\rangle/ \langle x_1+x_2+\cdots +x_{n}\rangle$$ can be identified with $\Q\langle Y_{n-1}\rangle$, where $y_i=x_i-x_{i+1}$ for $1\leq i \leq n-1$.

\subsection{Solomon's descent algebra of $S_{n}$}
\label{sec:so}

Surprisingly, the key to prove the general result comes from the theory of the descent algebra of the symmetric group which we will recall here. 
Consider the set  $I=\{1,2,\ldots,n-1\}$. The descent set of $\sigma\in {\Sn}_n$ is the set 
$Des(\sigma)=\{i\in I\,|\sigma(i)>\sigma(i+1)\}$. 
For $K\subseteq I$, define 
\begin{equation*}
D_K=\sum_{\sigma\in {\Sn}_n\atop Des(\sigma)=K} \sigma.
\end{equation*}
The Solomon's descent algebra $\Sigma ({\Sn}_n)$ is a subalgebra of the group algebra $\Q {\Sn}_n$ with basis $\{D_K|K\subseteq I\}$ 
\cite{S}.
For a standard tableau $t$ of shape $\lambda\vdash n$  define $$z_t=\sum_{\sigma\in {\Sn}_n\atop Q(\sigma)=t} \sigma,$$ where $Q(\sigma)$ corresponds to the recording tableau in the Robinson-Schensted correspondence.
Then consider the linear span ${\mathcal {Q}}_n={\mathcal{L}}\{z_t|\,\,t \in STab_n\}$. Note in general that ${\mathcal {Q}}_n$ is not a subalgebra of $\Q {\Sn}_n$, for $n\geq 4$.
Define the descent set of a standard tableau $t$ by
$Des(t)=\{i|i+1 {\text { is above i in }} t\}$.
We can observe that $Des(\sigma)=Des(\sigma')$ if and only if  $Des(Q(\sigma))=Des(Q(\sigma'))$ so we get the equality
 \begin{equation}\label{eq:zt}
 D_K=\sum_{t\in STab_n\atop Des(t)=K} z_t .
 \end{equation} 
 Hence $\Sigma({\Sn}_n)
\subseteq{\mathcal {Q}}_n$.
There is an algebra morphism from the descent algebra to the character algebra
$\theta:\Sigma ({\Sn}_n)\rightarrow \Q {\text{Irr}}({\Sn}_n)$
 due to Solomon \cite{S}.
Moreover, 
there is a linear map \cite{PR}
\begin{equation}\label{thetatilde}
\tilde{\theta}:{\mathcal {Q}}_n\rightarrow \Q {\text{Irr}}({\Sn}_n)
\end{equation}
defined by $\tilde{\theta}(z_t)=\chi^{{\rm shape}(t)}$, 
and $\tilde{\theta}$ restricted to $\Sigma({\Sn}_n)$ corresponds to $\theta$.
Note that this map is not an homomorphism with respect to the internal product on  $\Q {\text{Irr}}({\Sn}_n)$, in fact, ${\mathcal {Q}}_n$ does not have a corresponding product.

\subsection{General method for ${\Sn}_n$}
\label{sec:de}

 We are developing a combinatorial method to determine the multiplicity of $V^{\lambda}$ in $V^{\otimes d}$, when $V$ is any ${\Sn}_n$-module. For this purpose, we will use the algebra morphism $\theta:\Sigma({\Sn}_n)\rightarrow  \Q{\text{Irr}}({\Sn}_n) $ introduced in section \ref{sec:so}.
 The next proposition says that this multiplicity is given as some coefficients in $f^d$, when $f$ is an element of the descent algebra $\Sigma({\Sn}_n)$ which is sent to  the character $\chi^V$ of $V$.

\begin{prop}
\label{prop:genSn}
Let  $V$ be an ${\Sn}_n$-module such that $\theta(f)=\chi^V$, for some $f\in\Sigma({\Sn}_n)$.
For $\lambda\vdash n$, the multiplicity of  $V^{\lambda}$ in $V^{\otimes d}$ is  equal to
$$\sum_{t\in STab_n \atop shape(t)=\lambda}[z_t] {f}^d,$$ where $[z_t] {f}^d$ is the coefficient of $z_t$ in  ${f}^d$.
\end{prop}
\begin{proof}
By equation \eqref{eq:zt}, we can write ${f}^d=\sum_{\lambda\vdash n}\sum_{t\in STab_n \atop shape(t)=\lambda}c_tz_t.$
Applying the linear map ${\tilde{\theta}}$ of equation \eqref{thetatilde}, we get
\begin{equation*}
{\tilde{\theta}}(f^d)=\sum_{\lambda\vdash n}\sum_{t\in STab_n \atop shape(t)=\lambda}c_t\tilde{\theta}(z_t)=\sum_{\lambda\vdash n}\sum_{t\in STab_n \atop shape(t)=\lambda}c_t \chi^{\lambda}.\\
\end{equation*}
Now by equation  \eqref{thetatilde} and hypothesis, ${\tilde{\theta}}(f^d)={\theta}(f^d)={\theta}(f)^{d}=(\chi^{V})^{d}$,
so the coefficient of $\chi^{\lambda}$ in $(\chi^{V})^{d}$ is $$\sum_{t\in STab_n \atop shape(t)=\lambda}c_t=\sum_{t\in STab_n \atop shape(t)=\lambda}[z_t]f^d.$$
\QED\end{proof}

Although next theorem is an easy consequence of Lemma \ref{lemCoeff} and Proposition \ref{prop:genSn}, it provides us with an interesting interpretation for the multiplicity of $V^{\lambda}$ in the $d$-fold Kronecker product of a ${\Sn}_n$-module. This multiplicity is the weighted sum of words in a particular Cayley graph of ${\Sn}_n$ which reduce to elements $\sigma_t$, where $\sigma_t$ has recording tableau $t$ of shape $\lambda$ in the Robinson-Schensted correspondence. Recall that the \emph{support} of an element $f$ of the group algebra $\Q {\Sn}_n$ is defined by ${\text{supp}}(f)=\{\sigma\in {\Sn}_n|[\sigma]f\neq 0\}$, where $[\sigma]f$ is the coefficient of $\sigma$ in $f$.

\begin{thm} \label{thm:genSn}
Let  $V$ be an ${\Sn}_n$-module such that $\theta(f)=\chi^V$, for some $f\in\Sigma({\Sn}_n)$.
For $\lambda\vdash n$, the multiplicity of  $V^{\lambda}$ in $V^{\otimes d}$ is  
$$\sum_{t\in STab_n \atop shape(t)=\lambda}\sum_{w\in w(\sigma_t,d; \Gamma)}\omega(w),$$
where $\sigma_t\in {\Sn}_n$ has recording tableau $Q(\sigma_t)=t$, $\Gamma=\Gamma({\Sn}_n,{\text{supp}}(f))$ with associated weight $\omega(\sigma)=[\sigma](f)$ for each $\sigma\in{\text{supp}}(f)$ and $w(\sigma_t,d; \Gamma)$ is the set of words of length $d$ which reduce to $\sigma_t$ in $\Gamma$. In particular, the multiplicity of the trivial module is $$\sum_{w\in w(e,d; \Gamma)}\omega(w).$$
\end{thm}
\begin{proof}
From Proposition \ref{prop:genSn}, the multiplicity of  $V^{\lambda}$ in $V^{\otimes d}$ is  
$$\sum_{t\in STab_n \atop shape(t)=\lambda}[z_t] {f}^d.$$ 
Since by definition $\sigma\in{\text{supp}}(z_t)$ if and only if $\sigma$ has recording tableau $t$ in the Robinson-Schensted correspondence, the coefficient of $z_t$ in ${f}^d$ is also the coefficient of  $\sigma_t$ in ${f}^d$ with $Q(\sigma_t)=t$ and the result follows from Lemma \ref{lemCoeff}.
\QED\end{proof}

\subsection{Decomposition of $T(V^{(n-1,1)})$ and words}
\label{sec:deess}

Since we are particularly interested in the geometric ${\Sn}_n$-module, we make explicit the following two corollaries respectively of Proposition \ref{prop:genSn} and Theorem \ref{thm:genSn}
needed to draw a connection between the multiplicity of $V^{\lambda}$ in $(V^{(n-1,1)})^{\otimes d}$ and words of length $d$ in a particular Cayley graph of ${\Sn}_n$.  
The first corollary relates the multiplicity of  $V^{\lambda}$ in $(V^{(n-1,1)})^{\otimes d}$ to certain coefficients in the product  of the basis element ${D_{\{1\}}}^d$ of the descent algebra $\Sigma({\Sn}_n)$, where $D_{\{1\}}$ is the sum of all elements in ${\Sn}_n$ having descent set \{1\}. The second one relates the multiplicity of  $V^{\lambda}$ in $(V^{(n-1,1)})^{\otimes d}$ to words in the Cayley graph of ${\Sn}_n$ with generators $(12),(132),\ldots, (1\,n\,\cdots\, 432)$.

\begin{cor}\label{prop:coeff1}
Let $\lambda\vdash n$.
The multiplicity of  $V^{\lambda}$ in $(V^{(n-1,1)})^{\otimes d}$ is  
$$\sum_{t\in STab_n \atop shape(t)=\lambda}[z_t] {D_{\{1\}}}^d,$$
where $[z_t] {D_{\{1\}}}^d$ is the coefficient of $z_t$ in  ${D_{\{1\}}}^d$.
\end{cor}
\begin{proof}
Since the element of the linear span ${\mathcal {Q}}_n={\mathcal{L}}\{z_t|\,\,t \in STab_n\}$ of section \ref{sec:so}
\begin{align*} 
z_{  \newdimen\Squaresize \Squaresize=6pt 
{\,\Young
 {\Box{\scriptscriptstyle{2}}&\cr
  \Box{\scriptscriptstyle 1}&\Box{\scriptscriptstyle 3}&\Box{\scriptscriptstyle 4}&\Box{\scriptscriptstyle\,\cdots }&\Box{\scriptscriptstyle n}\cr}}}&=2\,1\,3\,4\,\cdots\,n+3\,1\,2\,4\,\cdots\,n+4\,1\,2\,3\,\cdots\,n+\cdots+n\,1\,2\,3\,\cdots\,n-1\\
&=(12)+(132)+(1432)+\cdots + (1\,n\,\cdots 432)\\
\label{lemzt0} 
&=D_{\{1\}},
\end{align*}
we have ${\theta}({D_{\{1\}}})=\chi^{(n-1,1)}$ and the results follows from Proposition \ref{prop:genSn}.
\QED\end{proof}

\begin{cor}\label{prop:Mult}
Let $\lambda\vdash n$. The multiplicity of  $V^{\lambda}$ in $(V^{(n-1,1)})^{\otimes d}$ is  equal to 
$$\sum_{t\in STab_n \atop shape(t)=\lambda}|w(\sigma_t,d; \Gamma)|,$$ where $\sigma_t\in {\Sn}_n$ has recording tableau $t$, $\Gamma=\Gamma({\Sn}_n,\{(12),(132),\ldots, (1\,n\,\cdots\, 432)\})$ and $w(\sigma_t,d; \Gamma)$ is the set of words of length $d$ which reduce to $\sigma_t$ in $\Gamma$. In particular, the multiplicity of the trivial module is equal to the number of words of length $d$ which reduce to the identity in $\Gamma$.
\end{cor}
\begin{proof}
As in the proof of Corollary \ref{prop:coeff1}, ${\theta}({D_{\{1\}}})=\chi^{(n-1,1)}$. Now  $$\sum_{t\in STab_n \atop shape(t)=\lambda}\sum_{w\in w(\sigma_t,d; \Gamma)}\omega(w)=\sum_{t\in STab_n \atop shape(t)=\lambda}\sum_{w\in w(\sigma_t,d; \Gamma)}1=\sum_{t\in STab_n \atop shape(t)=\lambda}|w(\sigma_t,d; \Gamma)|$$ and the result follows from Theorem \ref{thm:genSn}.
\QED\end{proof}


\begin{ex}\label{exMult}
Using Corollary \ref{prop:coeff1}, the ${\Sn}_3$-module $(V^{(2,1)})^{\otimes 4}$ decomposes into simple modules as $3\,V^{(3)}\oplus  5\,V^{(2,1)}\oplus 3\,V^{(1,1,1)}$ since 
\begin{align*}
{D_{\{1\}}}^4&=3\,D_{\emptyset}+3\,D_{\{2\}}+2\,D_{\{1\}}+3\,D_{\{1,2\}}\\
&=3\,z_{\newdimen\Squaresize \Squaresize=6pt 
\,\Young
 {\Box{\scriptstyle 1}&
  \Box{\scriptstyle 2}&
  \Box{\scriptstyle 3}\cr}}\,
  +3\,z_{\newdimen\Squaresize \Squaresize=6pt 
\,\Young
 {\Box{\scriptstyle 3}\cr
  \Box{\scriptstyle 1}&
  \Box{\scriptstyle 2}\cr}}\,
  +2\,z_{\newdimen\Squaresize \Squaresize=6pt 
\,\Young
 {\Box{\scriptstyle 2}\cr
  \Box{\scriptstyle 1}&
  \Box{\scriptstyle 3}\cr}}
  +3\,z_{\newdimen\Squaresize \Squaresize=6pt 
\,\Young
 {\Box{\scriptstyle 3}\cr
  \Box{\scriptstyle 2}\cr
  \Box{\scriptstyle 1}\cr}}.\\
\end{align*}
These multiplicities can also be computed using Corollary \ref{prop:Mult} in the following way. 
Consider the Cayley graph $\Gamma=\Gamma({\Sn}_3, \{(12), (132)\})$ of Figure \ref{graphe6}
 \begin{figure}[ht]
\begin{center}
\includegraphics[scale=0.8]{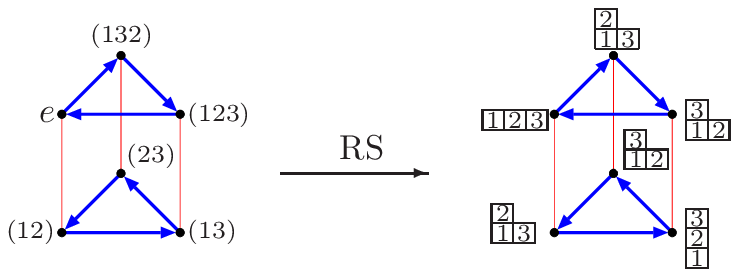}
\caption{{\bf Cayley graph $\Gamma({\Sn}_3, \{(12), (132)\})$ and recording tableaux}}
\label{graphe6}
\end{center}
\end{figure}
and write $a$ for $(12)$ and $b$ for $(132)$ to simplify. If we choose the representatives $$\sigma_{\newdimen\Squaresize \Squaresize=6pt 
\,\Young
 {\Box{\scriptstyle 1}&
  \Box{\scriptstyle 2}&
  \Box{\scriptstyle 3}\cr}}\,=e \qquad\sigma_{\newdimen\Squaresize \Squaresize=6pt 
\,\Young
 {\Box{\scriptstyle 3}\cr
  \Box{\scriptstyle 1}&
  \Box{\scriptstyle 2}\cr}}\,=(23) \qquad\sigma_{\newdimen\Squaresize \Squaresize=6pt 
\,\Young
 {\Box{\scriptstyle 2}\cr
  \Box{\scriptstyle 1}&
  \Box{\scriptstyle 3}\cr}}\,=(12)\qquad\sigma_{\newdimen\Squaresize \Squaresize=6pt 
\,\Young
 {\Box{\scriptstyle 1}\cr
  \Box{\scriptstyle 2}\cr
  \Box{\scriptstyle 3}\cr}}\,=(13)$$
the multiplicities of $V^{(3)}$, $V^{(2,1)}$ and $V^{(1,1,1)}$ are respectively  given by the cardinalities of  the sets of words (see Example \ref{ex:motsS31})
$$ \begin{array}{ll}
|w({e}\,,4;\Gamma)|  &= |\{aaaa, abab, baba\}|=3,\\
|w({(23)}\, ,4;\Gamma)|+|w({(12)}\,,4;\Gamma)| &=|\{aaba,baaa,bbab\}|+|\{abbb,bbba\}|=5,\\
|w({(13)}\,,4;\Gamma)|&= |\{aaab,abaa,babb\}|=3.
\end{array}$$
If instead we choose the representatives  
$$\sigma_{\newdimen\Squaresize \Squaresize=6pt 
\,\Young
 {\Box{\scriptstyle 3}\cr
  \Box{\scriptstyle 1}&
  \Box{\scriptstyle 2}\cr}}\,=(123)\qquad\sigma_{\newdimen\Squaresize \Squaresize=6pt 
\,\Young
 {\Box{\scriptstyle 2}\cr
  \Box{\scriptstyle 1}&
  \Box{\scriptstyle 3}\cr}}\,=(12)$$
 the mutiplicity of $V^{(2,1)}$ can also be computed by
 \begin{align*}
|w((123)\, ,4;\Gamma)|+| w((132)\,,4;\Gamma)|&=|\{aabb, baab, bbaa\}|+|\{abba, bbbb\}|=5.\\
  \end{align*}
These multiplicities can also be calculated using the master theorem for tensor products but it is not as combinatorial as by graphical means.  Using the inner product of characters and the character table for ${\Sn}_3$ (see Table \ref{tabcharS})
\begin{table}[ht!]
$$\begin{array}{|l||c|c|c|c|c|}
\hline
&(1,1,1)&(2,1)&(3)
\\\hline\hline
\chi^{(3)}&1&1&1\\\hline
\chi^{(2,1)}&2&0&-1\\\hline
\chi^{(1,1,1)}&1&-1&1\\\hline
\end{array}$$
\caption{\bf Character table of ${\Sn}_3$.}
\label{tabcharS}
 \end{table}
the multiplicities of $V^{(3)}$, $V^{(2,1)}$ and $V^{(1,1,1)}$ in $(V^{(2,1)})^{\otimes 4}$ are respectively computed by
\begin{align*}
\langle (\chi^{(2,1)})^4,\chi^{(3)} \rangle
&=\frac{1}{6}\bigg((\chi_{(1,1,1)}^{(2,1)})^4\chi_{(1,1,1)}^{(3)}+3\,(\chi_{(2,1)}^{(2,1)})^4\chi_{(2,1)}^{(3)})+2\,(\chi_{(3)}^{(2,1)})^4\chi_{(3)}^{(3)}\bigg)\\
&=\frac{1}{6}\big(2^4\cdot 1+3\cdot 0^4\cdot 1+2\cdot (-1)^4\cdot 1\big)=3,
\end{align*}
\begin{align*}
\langle (\chi^{(2,1)})^4,\chi^{(2,1)} \rangle
&=\frac{1}{6}\bigg((\chi_{(1,1,1)}^{(2,1)})^4\chi_{(1,1,1)}^{(2,1)}+3\,(\chi_{(2,1)}^{(2,1)})^4\chi_{(2,1)}^{(2,1)})+2\,(\chi_{(3)}^{(2,1)})^4\chi_{(3)}^{(2,1)}\bigg)\\
&=\frac{1}{6}\big( 2^4 \cdot 2+3\cdot 0^4\cdot 0+2\cdot (-1)^4\cdot (-1)\big)=5,
\end{align*}
\begin{align*}
\langle (\chi^{(2,1)})^4,\chi^{(1^3)} \rangle
&=\frac{1}{6}\bigg((\chi_{(1,1,1)}^{(2,1)})^4\chi_{(1,1,1)}^{(1,1,1)}+3\,(\chi_{(2,1)}^{(2,1)})^4\chi_{(2,1)}^{(1,1,1)})+2\,(\chi_{(3)}^{(2,1)})^4\chi_{(3)}^{(1,1,1)}\bigg)\\
&=\frac{1}{6}\big( 2^4 \cdot 1+3\cdot 0^4\cdot (-1)+2\cdot (-1)^4\cdot 1\big)=3.
\end{align*}
\end{ex}

\subsection{Invariant algebra $T(V^{(n-1,1)})^{{\Sn}_n}\simeq\Q\langle Y_{n-1}\rangle^{{\Sn}_n}$}\label{sec:invess}
In particular, we have an interpretation for the graded dimensions of the space $T(V^{(n-1,1)})^{{\Sn}_n}$ of invariants of ${\Sn}_n$ in terms of paths which begin and end at the identity vertex in the Cayley graph of ${\Sn}_n$ with generators $(12)$, $(132)$, $\ldots$,  $(1\,n\,\cdots\,432)$. As a corollary of Corollary \ref{prop:Mult}, we can show that the dimension of  $T(V^{(n-1,1)})^{{\Sn}_n}$ in each degree $d$ can be counted by those words of length $d$. 

\begin{cor}\label{cor:invess}
The dimension of $((V^{(n-1,1)})^{\otimes d})^{{\Sn}_n}\simeq\Q\langle Y_{n-1}\rangle_d^{{\Sn}_n}$ is equal to the number of words of length $d$ which reduce to the identity in the Cayley graph $\Gamma({\Sn}_n,\{(12), \,(132),\, \ldots, (1\,n\,\cdots\,432) \})$.
\end{cor}
\begin{proof} 
The dimension of the space of invariants in $\Q\langle Y_{n-1}\rangle_d\simeq(V^{(n-1,1)})^{\otimes d}$ is equal to the multiplicity of the trivial module in $(V^{(n-1,1)})^{\otimes d}$. Then the result follows from Corollary \ref{prop:Mult}.
  \QED\end{proof}

\begin{ex}
Using the Reynold's operator $\sum_{\sigma\in {\Sn}_n} \sigma$ acting on the monomials, a basis for the space $\Q\langle y_1,y_2 \rangle_4^{{\Sn}_3}$ of invariants of ${\Sn}_3$ is given by the three following polynomials
\begin{align*}
p_1(y_1,y_2)&=y_1^2y_2^2-y_1y_2^2y_1-y_2y_1^2y_2+y_2^2y_1^2,\\
p_2(y_1,y_2)&=y_1y_2y_1y_2-y_1y_2^2y_1-y_2y_1^2y_2+y_2y_1y_2y_1,\\
p_3(y_1,y_2)&=2y_1^4+y_1^3y_2+y_1^2y_2y_1+y_1y_2y_1^2+3y_1y_2^2y_1+y_1y_2^3\\&\qquad+y_2y_1^3+3y_2y_1^2y_2+y_2y_1y_2^2+y_2^2y_1y_2+y_2^3y_1+2y_2^4,
\end{align*}
which agree with the number of words $\{aaaa,abab,baba\}$ (see Example \ref{ex:motsS31}) in the letters $a=(12)$ and $b=(132)$ which reduce to the identity in the Cayley graph $\Gamma({\Sn}_3,\{(12), \,(132) \})$ of Figure \ref{graphe4}.
\end{ex}

Another interesting result is that the number of free generators of the algebra of  invariants of ${\Sn}_n$ can be counted by some special paths in the Cayley graph of ${\Sn}_n$. These are the paths which begin and end at the identity vertex, but without crossing the identity vertex. 
Since in general  the algebra $T(V)^G$ of invariants of $G$  is freely generated \cite{DF}, we have a relation between the Hilbert-Poincar\'e series $P(T(V)^G)$ and the generating series $F(T(V)^G)$ whose coefficients count the number of free generators  of $T(V)^G$, which is the following
\begin{equation}\label{equ:gen}
P(T(V)^G)=\frac{1}{1-F(T(V)^G)}.
\end{equation}

\begin{prop} 
The number of free generators of $T(V^{(n-1,1)})^{{\Sn}_n}$ as an algebra are counted by the words which reduce to the identity without crossing the identity in $\Gamma({\Sn}_n,\{(12), \,(132),\, \ldots, (1\,n\,\cdots\,432) \})$. \end{prop}
\begin{proof}
Let $\Gamma=\Gamma({\Sn}_n,\{(12), \,(132),\, \ldots, (1\,n\,\cdots\,432) \})$.
Let $A_{e}(q)$ be the generating series counting the number of words which reduce to the identity in $\Gamma$ and $B_{e}(q)$ the one counting the number of words which reduce to the identity  without crossing it in $\Gamma$.
From Corollary \ref{cor:invess}, we have $P(T(V^{(n-1,1)})^{{\Sn}_n})=A_{e}(q)$ and by Lemma \ref{lem:gen} and equation \eqref{equ:gen}, this equals  to
$$P(T(V^{(n-1,1)})^{{\Sn}_n})=A_{e}(q)=\frac{1}{1-B_e(q)}=\frac{1}{1-F(T(V^{(n-1,1)})^{{\Sn}_n})}.$$ So
the generating series giving the number of free generators of $T(V^{(n-1,1)})^{{\Sn}_n}$ is $$F(T(V^{(n-1,1)})^{{\Sn}_n})=B_e(q).$$\QED\end{proof}

\begin{ex}
The free generators of $T(V^{(2,1)})^{{\Sn}_3}$ are counted  by  the number of words  which reduce to the identity without crossing the identity in the Cayley graph $\Gamma({\Sn}_3,\{(12),(132)\})$ of Figure \ref{graphe4}. They are 
$$\{aa\}, \{bbb\},\{abab, baba\}, \{ abbba, baabb,bbaab\}, \ldots$$
with cardinalities corresponding to the Fibonacci numbers.
Indeed, using the analog of Molien's Theorem {\cite{DF}} and Lemma \ref{lem:gen},  
\begin{align*}
F(T(V^{(2,1)})^{{\Sn}_3})&=1-P(T(V^{(2,1)})^{{\Sn}_3})^{-1}\\
&=1-\Bigg(\frac{1}{6}\Big(\frac{1}{(1-2q)}+3+\frac{2}{(1+q)}\Big)\Bigg)^{-1}\\
&=1-\Bigg(\frac{1}{6}\Big(\frac{6-6q-6q^2}{1-q-2q^2}\Big)\Bigg)^{-1}\\
&=1-\Bigg(\frac{1-q-2q^2}{1-q-q^2}\Bigg)=\frac{q^2}{1-q-q^2},
\end{align*}
which is the generating series for the Fibonacci numbers.
\end{ex}

\subsection{Decomposition of $T(V^{(n)}\oplus V^{(n-1,1)})$ and words}
\label{sec:deness}

We were also interested by studying the case of  the permutation module $V^{(n)}\oplus V^{(n-1,1)}$.
Like in the case of  the geometric one, we can give an interpretation for the decomposition of $(V^{(n)}\oplus V^{(n-1,1)})^{\otimes d}$ in terms of words of length $d$ in the previous Cayley graph of ${\Sn}_n$, but with adding a loop corresponding to the identity to each vertex of the graph. Let us first make explicit the following corollary of Proposition \ref{prop:genSn} which relates the multiplicity of  $V^{\lambda}$ in $(V^{(n)}\oplus V^{(n-1,1)})^{\otimes d}$ to certain coefficients in the product  of the basis element $(e+{D_{\{1\}}})^d$ of the descent algebra $\Sigma({\Sn}_n)$.

\begin{cor}\label{prop:coeff2}
Let $\lambda\vdash n$.
The multiplicity of  $V^{\lambda}$ in $(V^{(n)}\oplus V^{(n-1,1)})^{\otimes d}$ is  
$$\sum_{t\in STab_n \atop shape(t)=\lambda}[z_t] {(e+D_{\{1\}})}^d,$$ where $[z_t] {(e+D_{\{1\}})}^d$ is the coefficient of $z_t$ in $ {(e+D_{\{1\}})}^d$.
\end{cor}
\begin{proof}
Since the element of the linear span ${\mathcal {Q}}_n={\mathcal{L}}\{z_t|\,\,t \in STab_n\}$ of section \ref{sec:so}
\begin{align*} 
z_{ \newdimen\Squaresize \Squaresize=6pt 
{\,\Young
 {
  \Box{\scriptscriptstyle 1}&\Box{\scriptscriptstyle 2}&\Box{\scriptscriptstyle 3}&\Box{\scriptscriptstyle\,\cdots }&\Box{\scriptscriptstyle n}\cr}}}
  &=1\,2\,3\,\cdots\,n=e,
\end{align*}
we have ${\theta}(e+{D_{\{1\}}})=\theta(e)+\theta(D_{\{1\}})=\chi^{(n)}+\chi^{(n-1,1)}$ and the result follows from Proposition \ref{prop:genSn}.
\QED\end{proof}

\begin{cor}\label{prop:Mult2} Let $\lambda\vdash n$. The multiplicity of $V^{\lambda}$ in $(V^{(n)}\oplus V^{(n-1,1)})^{\otimes d}$ is equal to 
  $$\sum_{t\in STab_n \atop shape(t)=\lambda}|w(\sigma_t,d;\Gamma)|,$$ where $\sigma_t\in {\Sn}_n$ has recording tableau $t$,  
$\Gamma=\Gamma({\Sn}_n,\{e,(12),(132),\ldots, (1\,n\,\cdots\, 432)\})$ and $w(\sigma_t,d;\Gamma)$ is the set of words of length $d$ which reduce to $\sigma_t$ in $\Gamma$. In particular, the multiplicity of the trivial module is $|w(e,d;\Gamma)|$.
\end{cor}
\begin{proof}
As in the proof of Corollary \ref{prop:coeff2}, ${\theta}(e+{D_{\{1\}}})=\chi^{(n)}+\chi^{(n-1,1)}$. The weight of each word is one, so the sum $$\sum_{w\in w(\sigma_t,d;\Gamma)}\omega(w)$$ is just the cardinality of  $w(\sigma_t,d;\Gamma)$ and the result follows from Theorem \ref{thm:genSn}.
\QED\end{proof}

\begin{ex} %
Using Corollary \ref{prop:coeff2}, the ${\Sn}_3$-module $(V^{(3)}\oplus V^{(2,1)})^{\otimes 3}$ decomposes into simple modules as $5\,V^{(3)}\oplus \, 9\,V^{(2,1)}\oplus \,4\,V^{(1,1,1)}$
 since
\begin{align*}
(D_{\emptyset}+{D_{\{1\}}})^3&=5\,D_{\emptyset}+5\,D_{\{2\}}+4\,D_{\{1\}}+4\,D_{\{1,2\}}\\
&=5\,z_{\newdimen\Squaresize \Squaresize=6pt 
\,\Young
 {\Box{\scriptstyle 1}&
  \Box{\scriptstyle 2}&
  \Box{\scriptstyle 3}\cr}}\,
  +5\,z_{\newdimen\Squaresize \Squaresize=6pt 
\,\Young
 {\Box{\scriptstyle 3}\cr
  \Box{\scriptstyle 1}&
  \Box{\scriptstyle 2}\cr}}\,
  +4\,z_{\newdimen\Squaresize \Squaresize=6pt 
\,\Young
 {\Box{\scriptstyle 2}\cr
  \Box{\scriptstyle 1}&
  \Box{\scriptstyle 3}\cr}}
  +4\,z_{\newdimen\Squaresize \Squaresize=6pt 
\,\Young
 {\Box{\scriptstyle 3}\cr
  \Box{\scriptstyle 2}\cr
  \Box{\scriptstyle 1}\cr}}.\\
\end{align*}
These multiplicities can also be computed using Corollary \ref{prop:Mult2} in the following way. Consider the Cayley graph $\Gamma=\Gamma({\Sn}_3,\{e,(12),(132)\})$ of Figure \ref{graphe7} 
\begin{figure}[ht] 
\begin{center}
\includegraphics[scale=0.8]{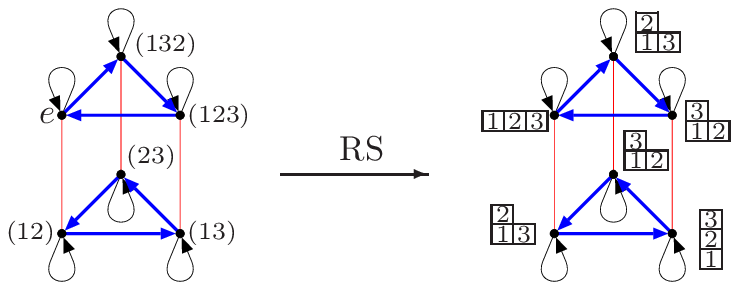}
\caption{{\bf Cayley graph $\Gamma({\Sn}_3,\{e,(12),(132)\})$ and recording tableaux.}}
\label{graphe7} 
\end{center}
\end{figure}
and set $a=(12)$ and $b=(132)$ to simplify.  If we choose the representatives $$\sigma_{\newdimen\Squaresize \Squaresize=6pt 
\,\Young
 {\Box{\scriptstyle 1}&
  \Box{\scriptstyle 2}&
  \Box{\scriptstyle 3}\cr}}\,=e \qquad\sigma_{\newdimen\Squaresize \Squaresize=6pt 
\,\Young
 {\Box{\scriptstyle 3}\cr
  \Box{\scriptstyle 1}&
  \Box{\scriptstyle 2}\cr}}\,=(123) \qquad\sigma_{\newdimen\Squaresize \Squaresize=6pt 
\,\Young
 {\Box{\scriptstyle 2}\cr
  \Box{\scriptstyle 1}&
  \Box{\scriptstyle 3}\cr}}\,=(132)\qquad\sigma_{\newdimen\Squaresize \Squaresize=6pt 
\,\Young
 {\Box{\scriptstyle 1}\cr
  \Box{\scriptstyle 2}\cr
  \Box{\scriptstyle 3}\cr}}\,=(13)$$
the multiplicities of $V^{(3)}$, $  V^{(2,1)}$ and $V^{(1,1,1)}$ are respectively  given by the cardinalities of  the following sets of words (see Example \ref{ex:motsS32})
\begin{align*}
|w(e,3;\Gamma)|
&= |\{bbb, aae, aea, eaa, eee\}|=5,\\
|w((123),3;\Gamma) |+\,|w((132),3;\Gamma)|&=|\{aba,bbe,beb,ebb\}|+ |\{aab,baa,bee,ebe,eeb\}|\\&=9,\\
|w((13),3;\Gamma)|&= |\{bba,abe,aeb,eab\}|=4.
\end{align*}
If instead we choose the representatives 
$$
\sigma_{\newdimen\Squaresize \Squaresize=6pt 
  \,\Young
 {\Box{\scriptstyle 3}\cr
  \Box{\scriptstyle 1}&
  \Box{\scriptstyle 2}\cr}}\,=(23) \qquad\sigma_{\newdimen\Squaresize \Squaresize=6pt 
\,\Young
 {\Box{\scriptstyle 2}\cr
  \Box{\scriptstyle 1}&
  \Box{\scriptstyle 3}\cr}}\,=(12)$$ 
 the mutiplicity of $V^{(2,1)}$ can also be computed by
\begin{align*}
|w((23),3;\Gamma)|+|w((12),3;\Gamma)|&=|\{abb,eab,bea,bae\}|+|\{aaa,bab,eea,eae,aee\}|\\&=9.
 \end{align*}
We will see later in section \ref{sec:ap} that the coefficient of the identity in $(e+{D_{\{1\}}})^d$ will also be the number of set partitions of $d$ into less than or equal to three parts because it is the dimension of the homogeneous symmetric polynomials of degree $d$ in three non-commutative variables.
\end{ex}

\subsection{Invariant algebra $T(V^{(n)}\oplus V^{(n-1,1)})^{{\Sn}_n}\simeq\Q\langle X_{n}\rangle^{{\Sn}_n}$}\label{sec:invness}

We give in the next corollary of Corollary \ref{prop:Mult2} the dimension of the homogeneous component of degree $d$ of $T(V)^{{\Sn}_n}\simeq\Q\langle X_n \rangle^{{\Sn}_n}$ in terms of paths in the Cayley graph of ${\Sn}_n$ which begin and end at the identity vertex. 

\begin{cor}
\label{cor:dim} 
The dimension of $((V^{(n)}\oplus V^{(n-1,1)})^{\otimes d})^{{\Sn}_n}\simeq\Q\langle X_{n}\rangle_d^{{\Sn}_n}$ is equal to the number of words of length $d$ which reduce to the identity in the Cayley graph $\Gamma({\Sn}_n,\{e, (12), (132),\, \ldots, (1\,n\,\cdots\,432) \})$.
\end{cor}
\begin{proof} 
The dimension of the space of invariants of ${\Sn}_n$ in $(V^{(n)}\oplus V^{(n-1,1)})^{\otimes d}\simeq\Q\langle X_{n}\rangle_d$ is equal to the multiplicity of the trivial module in $(V^{(n)}\oplus V^{(n-1,1)})^{\otimes d}$. Then the result follows from Corollary \ref{prop:Mult2}.
  \QED\end{proof}

\begin{ex}
A basis for the space $\Q\langle x_1,x_2,x_3 \rangle_3^{{\Sn}_3}$ is given by the five  following non-commutative monomial polynomials indexed by set partitions of $[3]$
\begin{align*}
&{\bf{m}}_{\{\{1\},\{2\},\{3\}\}}=x_1x_2x_3,\\
&{\bf{m}}_{\{\{1,2\},\{3\}\}}=x_1^2x_2+x_1^2x_3+x_2^2x_1+x_2^2x_3+x_3^2x_1+x_3^2x_2,\\
&{\bf{m}}_{\{\{1,3\},\{2\}\}}=x_1x_2x_1+x_1x_3x_1+x_2x_1x_2+x_2x_3x_2+x_3x_1x_3+x_3x_2x_3,\\
&{\bf{m}}_{\{\{1\}, \{2,3\}\}}=x_1x_2^2+x_1x_3^2+x_2x_1^2+x_2x_3^2+x_3x_1^2+x_3x_2^2,\\
&{\bf{m}}_{\{\{1,2,3\}\}}=x_1^3+x_2^3+x_3^3,\\
\end{align*}
\normalsize
which agree with the number of words $\{bbb,aae,aea,eaa,eee\}$ in the letters $a=(12)$, $b=(132)$ and $e$ which reduce to the identity in $\Gamma({\Sn}_3,\{e,(12), \,(132)\})$ of Figure \ref{graphe4}.
\end{ex}

Another interesting result is that the free generators of the invariant algebra are counted by some special paths in the Cayley graph of ${\Sn}_n$. These are the paths which begin and end at the identity vertex, but without crossing the identity vertex. 

\begin{prop} 
\label{prop:free}
The number of free generators of $T(V^{(n)}\oplus V^{(n-1,1)})^{{\Sn}_n}$ as an algebra are counted by the words  which reduce to the identity without crossing the identity in $\Gamma({\Sn}_n,\{e,(12), \,(132),\, \ldots, (1\,n\,\cdots\,432) \})$. 
\end{prop}
\begin{proof}
Follows from  Lemma \ref{lem:gen} and Equation \eqref{equ:gen}.
\QED\end{proof}

\begin{ex}
The free generators of $T(V^{(3)}\oplus V^{(2,1)})^{{\Sn}_3}$ are counted by the words which reduce to the identity without crossing the identity in the Cayley graph $\Gamma({\Sn}_3,\{e,(12),(132))\})$ of Figure \ref{graphe4}. They are
$$\{e\},\{aa\},\{bbb,aea\},\{abab,baba,bebb,bbeb,aeea\},\ldots$$
with cardinalities corresponding to the odd Fibonacci numbers beginning from the second subset.  Indeed, using the analog of Molien's Theorem and Lemma \ref{lem:gen}, 
\begin{align*}
F(T(V^{(3)}\oplus V^{(2,1)})^{{\Sn}_3})&=1-P(T(V^{(3)}\oplus V^{(2,1)})^{{\Sn}_3})^{-1}\\
&=1-\Bigg(\frac{1}{6}\Big(\frac{1}{(1-3q)}+2+\frac{3}{(1-q)}\Big)\Bigg)^{-1}\\
&=1-\Bigg(\frac{1}{6}\Big(\frac{6-18q+6q^2}{1-4q+3q^2}\Big)\Bigg)^{-1}\\
&=1-\Bigg(\frac{1-4q+3q^2}{1-3q+q^2}\Bigg)=\frac{q(1-2q)}{1-3q+q^2},
\end{align*}
which is the generating series for the odd Fibonacci numbers.
\end{ex}

\subsection{Applications to set partitions}
\label{sec:ap}
Let  $V=V^{(n)}\oplus V^{(n-1,1)}$ be the permutation module. The algebra  $T(V)^{{\Sn}_n}\simeq\Q\langle X_{n}\rangle^{{\Sn}_n}$ corresponds to the algebra of symmetric polynomials in non-commutative variables \cite{W,RS} with bases indexed by set partitions with at most $n$ parts.  That algebra is also freely generated by the set of non-commutative monomial polynomials indexed by non-splittable set partitions with at most $n$ parts \cite{BRRZ}.
Let $[n]=\{1,2,\ldots n\}$. A \emph {set partition of $[n]$}, denoted by $A\vdash [n]$, is a family of disjoint nonempty subsets $A_1, A_2, \ldots, A_k\subseteq [n]$ such that $A_1\cup A_2\cup\ldots\cup A_k=[n]$. The subsets $A_i$ are called the \emph{parts} of A.
Given two set partitions $B=\{B_1,B_2,\ldots B_k\}\vdash [n]$ and $C=\{C_1,C_2,\ldots  C_{\ell}\}\vdash [m]$, define
\begin{equation*}
B\circ C=
\begin{cases}
\{B_1\cup(C_1+n),\ldots,B_k\cup(C_k+n), (C_{k+1}+n),\ldots, (C_{\ell}+n)\} & \text{if $k\leq \ell$}\\
\{B_1\cup(C_1+n),\ldots,B_{\ell}\cup(C_{\ell}+n), B_{\ell+1}, \ldots, B_k\}&\text{if $k>\ell$}.\\
\end{cases}
\end{equation*}
A set partition $A$ is said to be {\emph{splitable}} if $A=B\circ C$, where $B$ and $C$ are non empty. 
\begin{ex}The set partitions
$\{\{1\},\{2\},\{3\}\}$ and $\{\{1\},\{2,3\}\}$ are nonsplitable and the remaining set partitions of $[3]$ split as
$$\begin{array}{l}
\{\{1,2\},\{3\}\}=\{\{1\}\}\circ\{\{1\},\{2\}\},\\
\{\{1,3\},\{2\}\}=\{\{1\},\{2\}\}\circ\{\{1\}\},\\
\{\{1,2,3\}\}=\{\{1\}\}\circ\{\{1\}\}\circ\{\{1\}\}.\\
\end{array}$$
\end{ex}

In Comtet \cite{C}, we can see that for a fixed $k\geq 0$, the ordinary and exponential generating functions giving the number of set partitions with $k$ parts are respectively  
\begin{align*}
\sum_{l\geq 0}S(l,k) q^l&=\frac {{q}^{k}}{ (1-0\,q)\left( 1-q \right)  \left( 1-2\,q \right)  \cdots \left( 1-k\,q \right)},\\
 \sum_{l\geq 0}S(l,k) \frac{q^l}{l!}&=\frac{1}{k!}(e^q-1)^k,
 \end{align*}
 where $S(l,k)$ is the Stirling numbers of second kind. 
Since the dimension of $(V^{\otimes d})^{{\Sn}_n}$ is given by the number of set partitions of $[d]$ wih at most $n$ parts,  the ordinary and exponential Hilbert-Poincar\'e series for  $T(V)^{{\Sn}_n}$ are repectively given by 
\begin{align*}
 P(T(V)^{{\Sn}_n})&=\displaystyle{\sum_{k=0}^{n} \frac {{q}^{k}}{ (1-0\,q)\left( 1-q \right)  \left( 1-2\,q
 \right)  \cdots \left( 1-k\,q \right) }},\\
 \tilde{P}(T(V)^{{\Sn}_n})&=\sum_{k=0}^{n}\frac{1}{k!}(e^q-1)^k.
\end{align*}
From a result of Chauve and Goupil \cite{CG}, the exponential Hilbert-Poincar\'e series
 of $T(V^{(n-1,1)})^{{\Sn}_n}$ is given by
\begin{equation} \tilde{P}(T(V^{(n-1,1)})^{{\Sn}_n})
=\sum_{k=0}^{n}\frac{(e^q-1)^k}{k!e^q}. \end{equation}
\noindent
Note that when $n$ goes to infinity, we have $$\sum_{k\geq0}\frac{(e^q-1)^k}{k!e^q}=\frac{e^{e^q-1}}{e^q}=e^{e^q-q-1},$$
which is equal to the generating series giving the number of set partitions without singleton. In other words, the coefficient of $q^d/d!$ in $e^{e^q-q-1}$ corresponds to the number of set partitions of $d$ into blocks of size greater that one.
\noindent We present here a conjecture for a closed formula giving the ordinary Hilbert-Poincar\'e series of $T(V^{(n-1,1)})^{{\Sn}_n}$ which does not seem to obviously follow from our combinatorial interpretations for the dimensions in terms of words in a Cayley graph of ${\Sn}_n$.
 \begin{conj} 
The ordinary Hilbert-Poincar\'e series of $T(V^{(n-1,1)})^{{\Sn}_n}$ is  
 \begin{align*}
 P(T(V^{(n-1,1)})^{{\Sn}_n})&=\frac{1}{1+q}+\frac{q}{1+q}\displaystyle{\sum_{k=0}^{n-1} \frac {{q}^{k}}{ (1-q)\left( 1-2\,q \right)  \cdots \left( 1-k\,q \right) }}\\
 &=\frac{1}{1+q}+\frac{q}{1+q}P(T(V^{(n-1)}\oplus V^{(n-2,1)})^{S_{n-1}}).\end{align*}
 \end{conj}
                 
            The interpretations in terms of words in a Cayley graph of ${\Sn}_n$ have a different nature to that of set partitions.  But from Corollary \ref{cor:dim} and Proposition \ref{prop:free} respectively, we can show for instance the following two corollaries.

\begin{cor}\label{cor:set}
The number of set partitions of $[d]$ into at most
$n$ parts is the number of words of length
$d$  which reduce to the identity  in the Cayley graph $\Gamma({\Sn}_n,\{e, (12), (132),\, \ldots, (1\,n\,\cdots\,432) \})$.
\end{cor}

\begin{cor}
The number of nonsplitable set partitions of $[d]$ into at most
$n$ parts is the number of words of length
$d$  in $\Gamma({\Sn}_n,\{e, (12), (132),\, \ldots, (1\,n\,\cdots\,432) \})$ which reduce to the identity without crossing the identity.
 \end{cor} 
               
\section{Dihedral group ${\D}_m$}
\label{sec:ot}

The same kind of results can be observed for other finite groups, for example in the case of cyclic and dihedral groups. We will present in this section the case of the dihedral group  ${\D}_m$  of order $2m$ with presentation ${\D}_m=\langle s,r\mid s^2 = r^m = srsr=e \rangle$. 
We will describe a combinatorial way to decompose the tensor algebra on any ${\D}_m$-module $V$ into simple modules by looking to words in a particular Cayley graph of ${\D}_m$. The bridge between those words and the decomposition of the tensor algebra $T(V)$ into simple modules is made possible via a subalgebra of the group algebra and a surjective morphism from that subalgebra into the algebra of characters. We will also make explicit some results about the algebra $T(V)^{{\D}_m}$ of invariants of ${\D}_m$ when $V$ is the geometric module.

\subsection{Simple ${\D}_m$-modules}
\label{sec:di}

For our purpose, let us first list the simple modules of the dihedral group ${\D}_m$.
\paragraph{For $m=2\,k$ even.} The conjugacy classes of ${\D}_m$ are $$\{e\}\quad\{r^{\pm 1}\}\quad\{r^{\pm 2}\}\quad\ldots\quad\{r^{\pm k}\}\quad\{r^{2i}s|0\leq i\leq k-1\}\quad\{r^{2i+1}s|0\leq i\leq k-1\}$$ hence up to isomorphisms, there are $k+3$ simple module $V^{id}$,  $V^{\gamma}$, $V^{\epsilon}$, $V^{\gamma\epsilon}$, and $V^{j}$,
for $1\leq j\leq k-1$, 
with associated characters $id$, $\gamma$, $\epsilon$, $\gamma\epsilon$ and $\chi_j$ (see Tables \ref{modulesDm} and \ref{charDm}).

\begin{table}[ht]
\begin{center}
\includegraphics[scale=0.8]{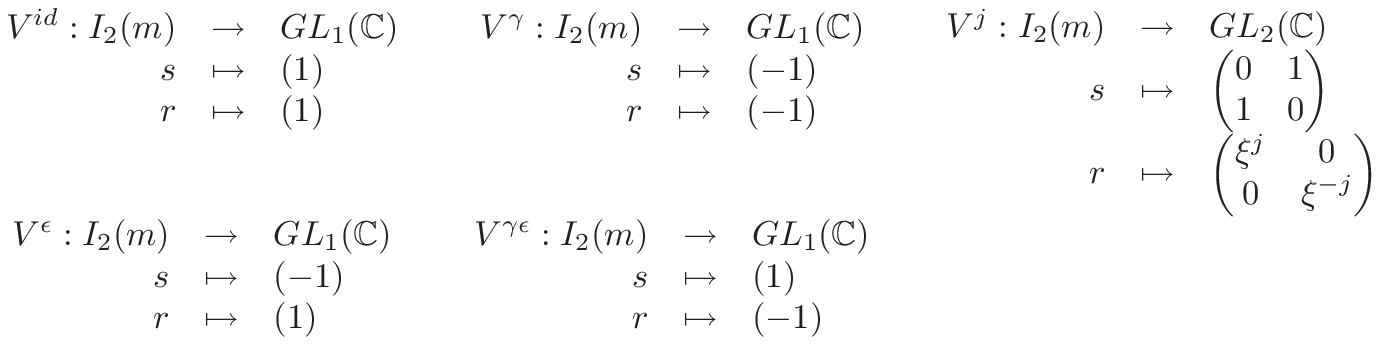}
\caption{{\bf Simple modules of ${\D}_m$, when $m$ is even and $\xi=e^{2\pi i/m}$.}}
\label{modulesDm}
\end{center}
\end{table}

\paragraph{For $m=2\,k+1$ odd.}  The conjugacy classes are $$\{e\}\quad\{r^{\pm 1}\} \quad\{r^{\pm 2}\} \quad\ldots \quad\{r^{\pm k}\} \quad\{r^{i}s|0\leq i\leq 2k\}$$ hence up to isomorphisms there are $k+2$ irreducible representations $V^{id}$, $V^{\epsilon}$ and $V^j$, for $1\leq j\leq k$, with associated irreducible characters  $id$, $\epsilon$ and  $\chi_j$ (see Tables \ref{modulesDm} and \ref{charDm}).

\begin{table}[ht]
\begin{center}
\includegraphics[scale=0.9]{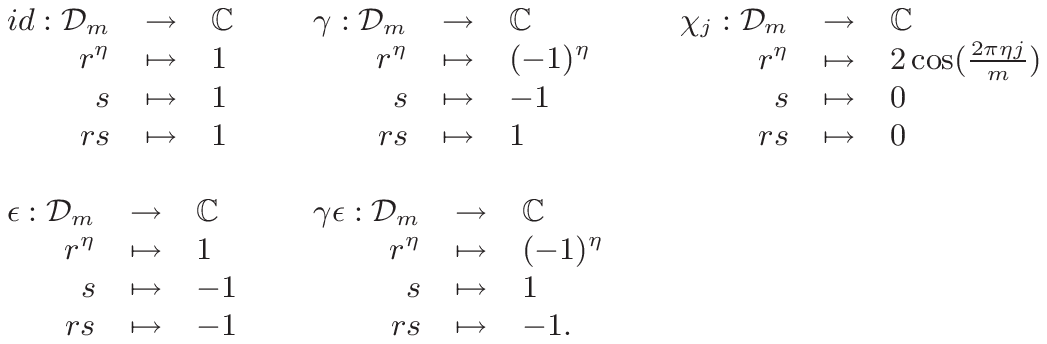}
\caption{{\bf Characters of ${\D}_m$, when $m$ is even.}}
\label{charDm}
\end{center}
\end{table}

\subsection{Subalgebra of the group algebra $\R {\D}_m$}
\label{sec:su}
The subalgebra of the group algebra that we will need is constructed in the next two propositions. 
We will consider two cases, when $m$ is even or odd, and will construct respectively two $\R$-linear spans with surjective algebra morphisms onto the algebra of characters.

\begin{prop} \label{prop:QD} Let $y_i=r^{1-i}s+r^{i}$. For $m=2k$ even, the linear span
 $${\mathcal{Q}}={\mathcal{L}}\{e,r^k,rs,r^{k+1}s, y_i, y_i rs\}_{1\leq i\leq k-1}$$ is a subalgebra of the group algebra $\R {\D}_m$, and there is a surjective algebra morphism 
$\theta:{\mathcal Q}\rightarrow \R {\rm Irr} ({\D}_m)$ defined by $\theta(e)=id$, $\theta(rs)=\epsilon$,  $\theta(r^k)= \gamma $, $\theta(r^{k+1}s)= \gamma\epsilon $ and $\theta(y_{i})=\theta( y_i rs)= \chi_i $.
\end{prop}

\begin{proof} 
The multiplication table of ${\mathcal Q}$ is represented by Table \ref{tabQeven}. Since $\chi_i\chi_j=\chi_j\chi_i$, we can assume $j\geq i$ and the multiplication table of  $\R \text{Irr}({\D}_m)$ is represented by Table \ref{tabChareven}.
\QED\end{proof}

 \begin{table}[ht!]
\begin{center}
\includegraphics[scale=0.8]{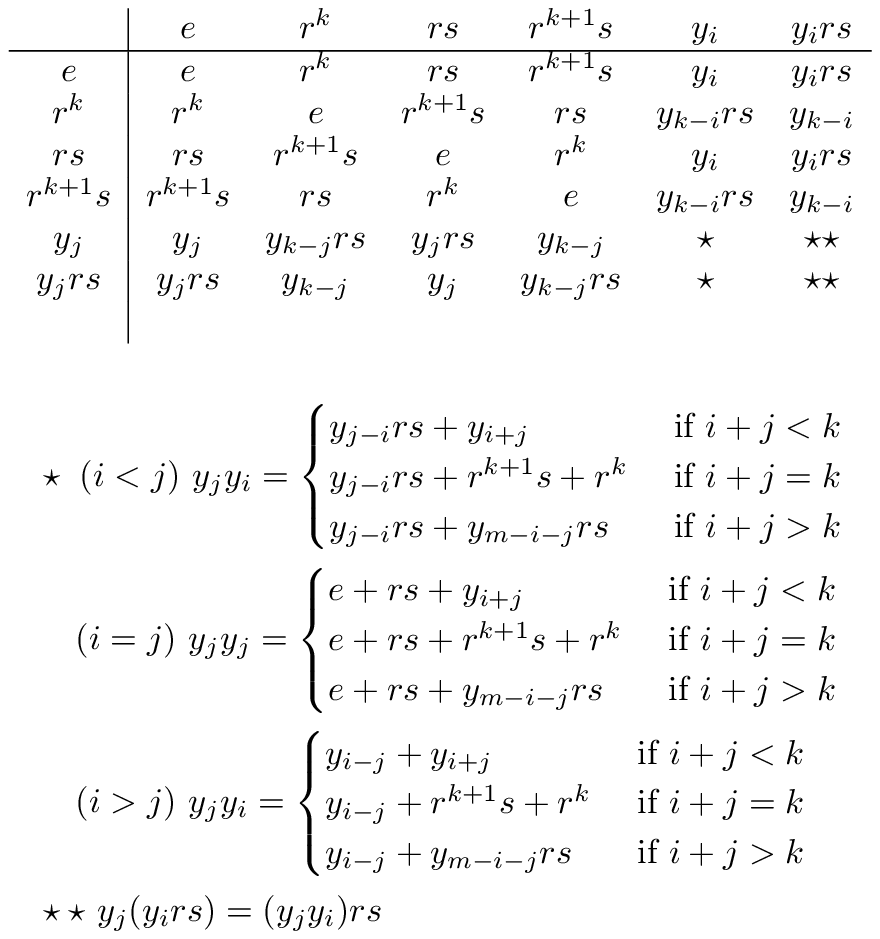}
   \caption{\bf Multiplication table of ${\mathcal{Q}}$, when $m=2k$ is even.}
   \label{tabQeven}
 \end{center}
 \end{table}
 
  \begin{table}[ht!]
\begin{center}
\includegraphics[scale=0.8]{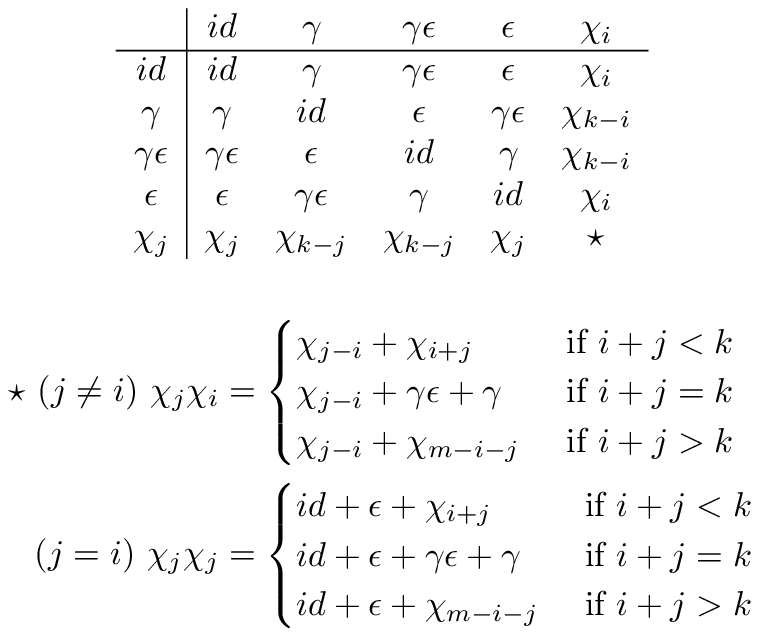}
   \caption{\bf Multiplication table of $\Z \text{Irr}({\D}_m)$, when $m=2k$ is even.}
    \label{tabChareven}
 \end{center}
 \end{table}

\begin{prop} Let $y_i=r^{1-i}s+r^{i}$. For $m=2k+1$ odd, the linear span
$${\mathcal{Q}}={\mathcal{L}}\{e, rs, y_i, y_i rs\}_{1\leq i\leq k}$$ is a subalgebra of the group algebra $\R {\D}_m$, and there is a surjective algebra morphism 
$\theta:{\mathcal Q}\rightarrow \R {\rm Irr} ({\D}_m)$ defined by 
$\theta(e)=id$, $\theta(rs)= \epsilon $ and  $\theta(y_{i})=\theta(y_i rs)= \chi_i$.
\end{prop}

\begin{proof}
See Table \ref{tabQodd} and Table \ref{tabCharodd}.
\QED\end{proof}

 \begin{table}[ht!]
\begin{center}
\includegraphics[scale=0.8]{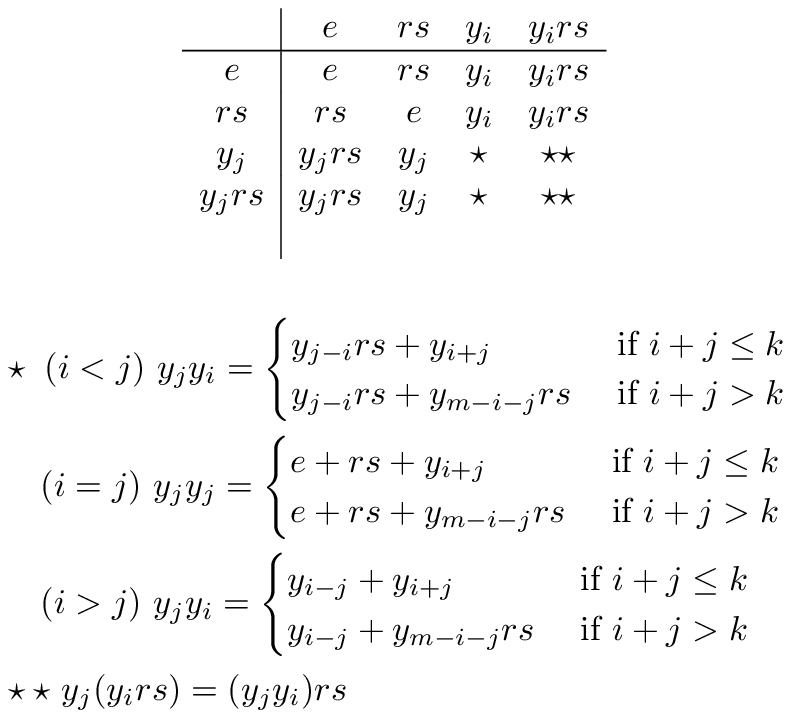}
   \caption{\bf Multiplication table of ${\mathcal{Q}}$, when $m=2k+1$ is odd.}
   \label{tabQodd}
 \end{center}
 \end{table}

\begin{table}[ht!]
\begin{center}
\includegraphics[scale=0.8]{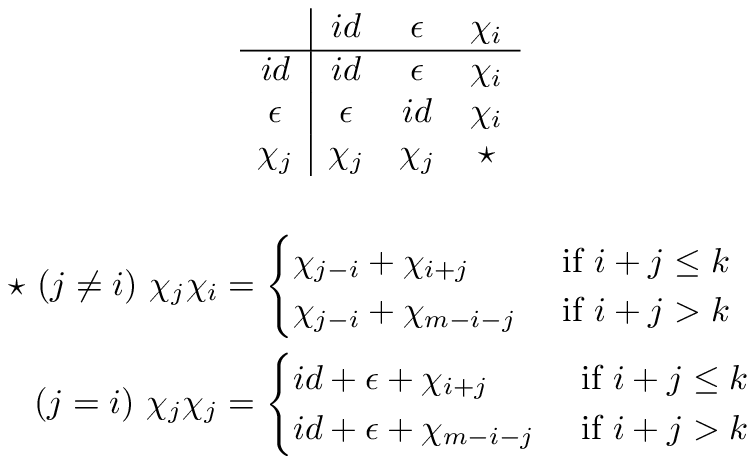}
   \caption{\bf Multiplication table of $\Z \text{Irr}({\D}_m)$, when $m=2k+1$ is odd.}
     \label{tabCharodd}
 \end{center}
 \end{table}

\subsection{General method for ${\D}_m$}
\label{sec:genD}

In this section is presented the general result for the dihedral groups.
To simplify the notation, $V^{(i)}$ will denote a simple ${\D}_m$-module with character $\chi^{(i)}$ and we will write the subalgebra of section \ref{sec:su} as ${\mathcal{Q}}={\mathcal{L}}\{b_i\}_{i\in I}$, where each element of the basis $b_i$ is sent to an irreducible character $\chi^{(i)}$ by the algebra morphism $\theta$. Recall also that the \emph{support} of an element $f$ of the group algebra $\R {\D}_m$ is defined by ${\text{supp}}(f)=\{g\in {\D}_m|[g]f\neq 0\}$, where $[g]f$ means taking the coefficient of $g$ in $f$. The first proposition relates the multiplicity of $V^{(k)}$ in $V^{\otimes d}$ to some coefficient in the expansion of the element $f^d$ of the subalgebra ${\mathcal{Q}}$ such that $\theta(f)=\chi^V$ and the next theorem  to words in a particuler Cayley graph of ${\D}_m$.

\begin{prop} \label{prop:coeffD}
Let $V$ be a ${\D}_m$-module.  If $f\in{\mathcal{Q}}$ is such that $\theta(f)=\chi^V$, then
the multiplicity of $V^{(k)}$ in $V^{\otimes d}$ is equal to 
$$\sum_{b_i \atop \theta(b_i)=\chi^{(k)}} [b_i]f^d,$$ where  $[b_i]f^d$ is the coefficient of  $b_i$ in $f^d$.
\end{prop}

\begin{proof} Let ${\mathcal{B}}=\{b_i\}_{i\in I}$.
Since $f^d$ is an element of ${\mathcal{Q}}$, we can write ${f}^d=\sum_{b_i\in {\mathcal{B}}}c_i b_i.$
Applying ${{\theta}}$ we get
\begin{equation*}
{{\theta}}(f^d)=\sum_{b_i\in {\mathcal{B}}}c_i {\theta(b_i)}=\sum_{b_i\in {\mathcal{B}}}c_i \chi^{(i)}.\\
\end{equation*}
On the left hand side, ${\theta}(f^d)={\theta}(f)^{d}=(\chi^{V})^{d}$,
so the coefficient of $\chi^{(k)}$ in $(\chi^{V})^{d}$ is equal to the sum of the coefficients of $b_i$ in $f^d$ such that $\theta(b_i)=\chi^{(k)}$.
\QED\end{proof}

\begin{thm}\label{thm:motsD} Let $V$ be a ${\D}_m$-module. If $f\in{\mathcal{Q}}$ is such that $\theta(f)=\chi^V$, then
the multiplicity of $V^{(k)}$ in $V^{\otimes d}$ is equal to 
$$\sum_{b_i \atop \theta(b_i)=\chi^{(k)}} \sum_{w\in w(\sigma_i,d;\Gamma)} \omega(w),$$
where $\sigma_i\in supp(b_i)$, $\Gamma=\Gamma({\D}_m, supp(f))$ with $\omega(g)=[g](f)$ for each $g\in supp(f)$ and $w(\sigma_i,d;\Gamma)$ is the set of words of length $d$ which reduce to $\sigma_i$ in  
$\Gamma$.
\end{thm}
\begin{proof}
Let ${\mathcal{B}}=\{b_i\}_{i\in I}$.
From Proposition \ref{prop:coeffD}, the multiplicity of $V^{(k)}$ in $V^{\otimes d}$ is equal to 
$$\sum_{b_i \in  {\mathcal{B}} \atop {\theta}(b_i)=\chi^{(k)}} [b_i] f^d.$$
Since each element in $supp(b_i)$ has coefficient one, we get $[b_i]f^d=[\sigma_i]f^d$. Moreover, all supports of the $b_i$'s are disjoint so using Lemma \ref{lemCoeff}  we get
$$[\sigma_i]f^d=\sum_{w\in w(\sigma_i,d;\Gamma)} \omega(w).$$ 
\QED\end{proof}

\begin{ex} Consider the ${\D}_4$-module $(2\,V^{1}\oplus V^{\gamma\epsilon})^{\otimes 2}$.  By Proposition  \ref{prop:QD}, there is a subagebra ${\mathcal{Q}}={\mathcal{L}}\{e,r^2,rs,r^{3}s, s+r, r^{3}+r^{2}s\}$ of $\R {\D}_4$ and an algebra morphism $\theta: {\mathcal Q}\rightarrow \R {\rm Irr} ({\D}_4)$ defined by 
$$\begin{array}{lll}
\theta(e)=id&\theta(r^2)= \gamma& \theta(s+r)= \chi_1\\
\theta(rs)=\epsilon&\theta(r^{3}s)= \gamma\epsilon&\theta( r^3+r^2s)= \chi_1.\\
\end{array}$$
Let $f=2\,( r^{3}+r^{2}s)+r^3s$. Applying $\theta$, the element
$f^2=5\,e+4\,rs+4\,r^2+4\,r^3s+2\,(s+r)+2\,( r^{3}+r^{2}s)$
is sent to 
$$(2\,\chi_1+\gamma\epsilon)^2={{5}}\,id+{{4}}\,\epsilon+{{4}}\,\gamma+{{4}}\,\gamma\epsilon+{{2}}\,\chi_1+{{2}}\,\chi_1$$ 
so the decomposition into simple modules is
$$(2\,V^{1}\oplus V^{\gamma\epsilon})^{\otimes 2}= {{5}}\,V^{id} \oplus {{4}}\,V^{\epsilon}\oplus {{4}}\,V^{\gamma} \oplus {{4}}\,V^{\gamma\epsilon}\oplus {{4}}V^{1} .$$
Using Theorem \ref{thm:motsD}, these multiplicities can be computed by counting weighted words in the Cayley graph
$\Gamma=\Gamma({\D}_4,\{{ {r^2s}}, { {r^3}},r^3s\})$ with weights $\omega(r^2s)=\omega(r^3)=2$ and $\omega(r^3s)=1$ of Figure \ref{graphe5}.
%
%
Set $a=r^2s$, $b=r^3$ and $c=r^3s$ to simplify. The multiplicities of $V^{id}$,  $V^{\epsilon}$,  $V^{\gamma}$,  $V^{\gamma\epsilon}$ and  $V^{1}$ are respectively given by 
$$ \begin{array}{ll}
V^{id}:&\displaystyle{\sum_{w\in w(e,2;\Gamma)}\omega(w)}=\omega(aa)+\omega(cc)=2\cdot 2+1\cdot 1=5\\ 
V^{\epsilon}:&\displaystyle{\sum_{w\in w(rs,2; \Gamma)}\omega(w)}=\omega(ba)=2\cdot 2=4\\
V^{\gamma}:&\displaystyle{\sum_{w\in w(r^2,2; \Gamma)}\omega(w)}=\omega(bb)=2\cdot 2=4\\
V^{\gamma\epsilon}:&\displaystyle{\sum_{w\in w(r^3s,2; \Gamma)}\omega(w)}=\omega(ab)=2\cdot 2=4\\
V^{1}:&\displaystyle{\sum_{w\in w(r,2; \Gamma)}\omega(w)}+\displaystyle{\sum_{w\in w(r^3,2; \Gamma)}\omega(w)}=\omega(ca)+\omega(ac)=1\cdot 2+2\cdot 1=4.\\
\end{array}$$
Note that  the multiplicity of $V^1$ can also be computed by
$$\displaystyle{\sum_{w\in w(s,2; \Gamma)}\omega(w)}+\displaystyle{\sum_{w\in w(r^2s,2; \Gamma)}\omega(w)}=\omega(cb)+\omega(bc)=1\cdot 2+2\cdot 1=4.$$
\end{ex}

\subsection{Invariant algebra $T(V^1)^{{\D}_m}$}
\label{sec:InvD}

We were particularly interested in studying the space of invariants of ${\D}_m$ in the tensor algebra on the geometric module $V^1$ and  we have the following results. The first one, which gives the graded dimensions of $T(V^1)^{{\D}_m}$ in terms of words, is a corollary of Theorem \ref{thm:motsD}.

\begin{cor}
The dimension of $((V^1)^{\otimes d})^{{\D}_m}\simeq\R\langle x_1,x_2\rangle_d^{{\D}_m}$ is equal to the number of words  of length $d$ wich reduce to the identity in the Cayley graph $\Gamma({\D}_m,\{r,s\})$.
\end{cor}
\begin{proof}
 The dimension of $\R\langle x_1,x_2\rangle_d^{{\D}_m}$ is equal to the multiplicity of the trivial module in $(V^1)^{\otimes d}$ and the result follows from Theorem \ref{thm:motsD} since $\theta(s+r)=\chi_1$.
 \QED\end{proof}

\begin{ex}
Consider the dihedral group ${\D}_4$ acting on $\R\langle x_1,x_2\rangle$ as
$$\begin{array}{ll}
s\cdot x_1=-x_1 & r \cdot x_1=x_1+\sqrt{2}\,x_2\\
s \cdot x_2=\sqrt{2}\, x_1+x_2& r \cdot x_2=-\sqrt{2}\,x_1-x_2
\end{array}.$$
Note that this action corresponds to the geometric action of ${\D}_4$.
Using Reynold's operator, a basis for $\R\langle x_1,x_2\rangle_4^{{\D}_4}$ is given by the four following polynomials
\begin{align*}
p_1(x_1,x_2)&=x_1x_2^2x_1+\frac{x_1x_2^3}{\sqrt{2}}+x_2x_1^2x_1+\frac{x_2x_1x_2^2}{\sqrt{2}}+\frac{x_2^2x_1x_2}{\sqrt{2}}+\frac{x_2^3x_1}{\sqrt{2}}+x_2^4,\\
p_2(x_1,x_2)&=x_1^2x_2^2+\frac{x_1x_2^3}{\sqrt{2}}+\frac{x_2x_1x_2^2}{\sqrt{2}}+x_2^2x_1^2+\frac{x_2^2x_1x_2}{\sqrt{2}}+\frac{x_2^3x_1}{\sqrt{2}}+x_2^4,\\
p_3(x_1,x_2)&=x_1x_2x_1x_2+\frac{x_1x_2^3}{\sqrt{2}}+x_2x_1x_2x_1+\frac{x_2x_1x_2^2}{\sqrt{2}}+\frac{x_2^2x_1x_2}{\sqrt{2}}+\frac{x_2^3x_1}{\sqrt{2}}+x_2^4,\\
p_4(x_1,x_2)&=x_1^4+\frac{x_1^3x_2}{\sqrt{2}}+\frac{x_1^2x_2x_1}{\sqrt{2}}+\frac{x_1x_2x_1^2}{\sqrt{2}}-\frac{x_1x_2^3}{\sqrt{2}}+\frac{x_2x_1^3}{\sqrt{2}}\\
&\qquad \qquad -\frac{x_2x_1x_2^2}{\sqrt{2}}-\frac{x_2^2x_1x_2}{\sqrt{2}}-\frac{x_2^3x_1}{\sqrt{2}}-x_2^4,\\
\end{align*}
\normalsize
which agree with the number of words of length four which reduce to the identity $\{ssss,rsrs,srsr,rrrr\}$ in the Cayley graph $\Gamma({\D}_4, \{r,s\})$ of Figure \ref{graphe9}.
\begin{figure}[ht]
\begin{center}
\includegraphics[scale=0.8]{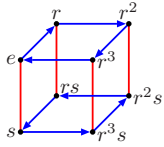}
\caption{{\bf Cayley graph $\Gamma({\D}_4, \{r,s\})$.}}
\label{graphe9}
\end{center}
\end{figure}

\end{ex}

As for the case of the symmetric group, we also have an interpretation for the free generators of $T(V^1)^{{\D}_m}$ as an algebra in terms of words.

\begin{prop} 
The number of free generators of $T(V^1)^{{\D}_m}$ as an algebra are counted by the words in the Cayley graph $\Gamma({\D}_m,\{r,s\})$ which reduce to the identity without crossing the identity. 
\end{prop}
\begin{proof}
Follows from  Lemma \ref{lem:gen} and Equation \eqref{equ:gen}.
\QED\end{proof}

Finally  we give a closed formula for the Hilbert-Poincar\'e series of the space $T(V^1)^{{\D}_m}$ of invariants of ${\D}_m$ in the next proposition.

\begin{prop} \label{prop: SerieI2}
The Hilbert-Poincar\'e series of $T(V^1)^{{\D}_m}\simeq\R\langle x_1,x_2\rangle^{{\D}_m}$  is 
 $$P(T(V^1)^{{\D}_m})=1+\frac{1}{2}
             \Bigg(\frac{(2q)^m+\sum_{i=0}^{\lfloor m/2\rfloor}(\binom{m+1}{2i+1}-2\binom{m}{2i}) (1-4q^2)^{i}}
             {\sum_{i=0}^{\lfloor m/2\rfloor}\binom{m}{2i}(1-4q^2)^i-(2q)^m}\Bigg).$$
\end{prop}

\begin{proof}
The Cayley graphs of ${\D}_3$, ${\D}_4$ and more generally ${\D}_m$ with generators $s$ and $r$ are represented in Figure \ref{graphe10}. When $m$ goes to infinity, we can visually represent it by the infinite strip of Figure \ref{graphe11}.
\begin{figure}[ht!]
\begin{center}
\includegraphics[scale=0.8]{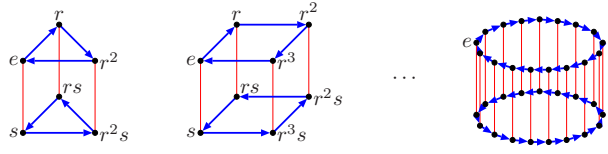}
\caption{{\bf Cayley graphs $\Gamma({\D}_3,\{s,r\})$, $\Gamma({\D}_4,\{s,r\})$ and more generally $\Gamma({\D}_m,\{s,r\})$.}}
\label{graphe10}
\end{center}
\end{figure}
\begin{figure}[ht!]
\begin{center}
\includegraphics[scale=0.8]{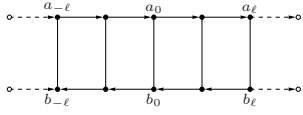}
\caption{{\bf Cayley graph $\Gamma({\D}_m,\{s,r\})$ when $m$ goes to infinity.}}
\label{graphe11}
\end{center}
\end{figure}

Then the words which reduce to the identity can be seen as the paths that go from $a_0$ to $a_0$, in addition to the paths from $a_0$ to $a_{\ell m}$ and  from $a_0$ to $a_{-\ell m}$, for $\ell\geq 0 $.
Let us define the generating functions of paths ending on the $a_i$'s side and the one of those ending on the $b_i$'s side as
\begin{align*}
A&=1+t\,qA+qB\\
B&=B/t+q A
\end{align*}
where $t$ counts the level from the starting point and $q$ the length.
Solving these equations we get
\begin{align*}
A(t,q)&=\sum_{\ell\in\Z,k\geq 0}c_{k,\ell}t^{\ell}\,q^k=\frac{q-t}{qt^2+q-t}\\
B(t,q)&=\sum_{\ell\in\Z,k\geq 0}d_{k,\ell}t^{\ell}\,q^k=\frac{-qt}{qt^2+q-t}\\
\end{align*}
where $c_{k,\ell}$ is the number of paths of length $k$ ending at $a_{\ell}$ and $d_{k,\ell}$ is the number of paths of length $k$ ending at $b_{\ell}$. To get the generating function counting the paths from $a_0$ to $a_{\ell}$ of length $k$, we need to take the coefficient of $t^{\ell}$ in $A(t,q)$
\begin{equation*}
[t^{\ell}] A(t,q)
=\sum_{k\geq 0} \binom{2k+\ell-1}{k} q^{2k+\ell}.
\end{equation*}
Therefore the number of paths from $a_0$ to $a_0$ is $\sum_{k\geq 0} \binom{2k-1}{k} q^{2k}$, the number of paths from $a_0$ to $a_{\ell m}$ is $\sum_{k\geq 0} \binom{2k+\ell m-1}{k} q^{2k+ \ell m}$ and the number of paths from $a_0$ to $a_{-\ell m}$ is $\sum_{k\geq 0} \binom{2k+\ell m-1}{k-1} q^{2k+\ell m}$. 
We then get
\begin{align}\label{eq1}
 P(T(V^1)^{{\D}_m})
 \nonumber&=\sum_{\ell\geq 1} \sum_{k\geq 0}\binom{2k+\ell m}{k}q^{2k+ \ell m}+\sum_{k\geq 0}\binom{2k-1}{k}q^{2k}\\
   \nonumber        &=\sum_{\ell\geq 1}\frac{(2q)^{ \ell m}}{(1+\sqrt{1-4q^2})^{\ell m}\sqrt{1-4q^2}}+\frac{1+\sqrt{1-4q^2}}{2\sqrt{1-4q^2}}\\
           &=\frac{1}{\sqrt{1-4q^2}}\Bigg(\sum_{\ell\geq 1}\frac{(2q)^{ \ell m}}{(1+\sqrt{1-4q^2})^{\ell m}}+\frac{1+\sqrt{1-4q^2}}{2}\Bigg).
\end{align}         
Since $GP=\sum_{\ell\geq 1}\frac{(2q)^{ \ell m}}{(1+\sqrt{1-4q^2})^{\ell m}}$ is a geometric progression we have the equalities
\begin{align}\label{eq2}
GP
           \nonumber&=\frac{(2q)^{m}}{(1+\sqrt{1-4q^2})^{m}-(2q)^m} \\
           \nonumber&=\frac{(2q)^{m}((1-\sqrt{1-4q^2})^{m}-(2q)^m)}{((1+\sqrt{1-4q^2})^{m}-(2q)^m)((1-\sqrt{1-4q^2})^{m}-(2q)^m)}\\
           &=\frac{(1-\sqrt{1-4q^2})^{m}-(2q)^m}{2(2q)^m-(1+\sqrt{1-4q^2})^{m}-(1-\sqrt{1-4q^2})^{m}}.
\end{align}
We also have the two identities
\begin{align}\label{eq3}
 (1+\sqrt{1-4q^2})^m \nonumber&=\sum_{i=0}^m \binom{m}{i}(\sqrt{1-4q^2})^i\\
\nonumber&=\sum_{i=0}^{\lfloor m/2\rfloor} \binom{m}{2i}(\sqrt{1-4q^2})^{2i}\\
&\qquad+\sum_{i=0}^{\lfloor (m-1)/2\rfloor} \binom{m}{2i+1}(\sqrt{1-4q^2})^{2i+1}
\end{align}
\begin{align}\label{eq4}
(1-\sqrt{1-4q^2})^m \nonumber&=\sum_{i=0}^m (-1)^i\binom{m}{i}(\sqrt{1-4q^2})^i\\
\nonumber&=\sum_{i=0}^{\lfloor m/2\rfloor} \binom{m}{2i}(\sqrt{1-4q^2})^{2i}\\
& \qquad-\sum_{i=0}^{\lfloor (m-1)/2\rfloor} \binom{m}{2i+1}(\sqrt{1-4q^2})^{2i+1}.
\end{align}
Therefore replacing \eqref{eq3} and  \eqref{eq4} in \eqref{eq2} gives
 \begin{align}         \label{eq5}
            GP\nonumber&=\frac{\sum_{i=0}^{\lfloor m/2\rfloor} \binom{m}{2i}(1-4q^2)^{i}-\sum_{i=0}^{\lfloor (m-1)/2\rfloor} \binom{m}{2i+1}(\sqrt{1-4q^2})^{2i+1}-(2q)^m}{2(2q)^m-2\sum_{i=0}^{\lfloor m/2\rfloor}\binom{m}{2i}(1-4q^2)^i} \\
         &=\frac{-\sum_{i=0}^{\lfloor (m-1)/2\rfloor} \binom{m}{2i+1}(\sqrt{1-4q^2})^{2i+1}}
             {2((2q)^m-\sum_{i=0}^{\lfloor m/2\rfloor}\binom{m}{2i}(1-4q^2)^i)}-\frac{1}{2}.
     \end{align}
Now substituting \eqref{eq5} in \eqref{eq1} we get
         \begin{align*}          
            P(T(V^1)^{{\D}_m})         
          &=\frac{1}{\sqrt{1-4q^2}}
             \Bigg(\frac{-\sum_{i=0}^{\lfloor (m-1)/2\rfloor} \binom{m}{2i+1}(\sqrt{1-4q^2})^{2i+1}}
             {2\big((2q)^m-\sum_{i=0}^{\lfloor m/2\rfloor}\binom{m}{2i}(1-4q^2)^i\big)}+\frac{\sqrt{1-4q^2}}{2}\Bigg) \\  
             &=\frac{-\sum_{i=0}^{\lfloor (m-1)/2\rfloor} \binom{m}{2i+1}(1-4q^2)^{i}+(2q)^m-\sum_{i=0}^{\lfloor m/2\rfloor}\binom{m}{2i}(1-4q^2)^i}
             {2\big((2q)^m-\sum_{i=0}^{\lfloor m/2\rfloor}\binom{m}{2i}(1-4q^2)^i\big)} \\      
            \end{align*}          
and using the binomial identity $\binom{n}{k}=\binom{n-1}{k-1}+\binom{n-1}{k} $, we finally get the desired result
\begin{align*}
P(T(V^1)^{{\D}_m})
             &=\frac{(2q)^m-\sum_{i=0}^{\lfloor (m-1)/2\rfloor} \binom{m+1}{2i+1}(1-4q^2)^{i}-\chi(m\,\,{\rm{even}})(1-4q^2)^{m/2}}
             {2\big((2q)^m-\sum_{i=0}^{\lfloor m/2\rfloor}\binom{m}{2i}(1-4q^2)^i\big)}\\
             &=\frac{1}{2}
             \Bigg(\frac{(2q)^m-\sum_{i=0}^{\lfloor m/2\rfloor} \binom{m+1}{2i+1}(1-4q^2)^{i}}
             {(2q)^m-\sum_{i=0}^{\lfloor m/2\rfloor}\binom{m}{2i}(1-4q^2)^i}\Bigg)\\  
             &=1+\frac{1}{2}
             \Bigg(\frac{(2q)^m+\sum_{i=0}^{\lfloor m/2\rfloor}(\binom{m+1}{2i+1}-2\binom{m}{2i}) (1-4q^2)^{i}}
             {\sum_{i=0}^{\lfloor m/2\rfloor}\binom{m}{2i}(1-4q^2)^i-(2q)^m}\Bigg). 
\end{align*}
\QED\end{proof}


\section{Appendix}
\label{ap}
The first four tables of this section represent the dimensions of the algebra $T(V)^{{\Sn}_n}$ of invariants  of the symmetric group, for $V$ the geometric module and the permutation one.  Similarly, the last two tables record the dimensions of the algebra $T(V)^{{\D}_m}$ of invariants of the dihedral group when $V$ is the geometric module.%

\begin{table}[ht!]
\begin{center} 
\includegraphics[scale=0.75]{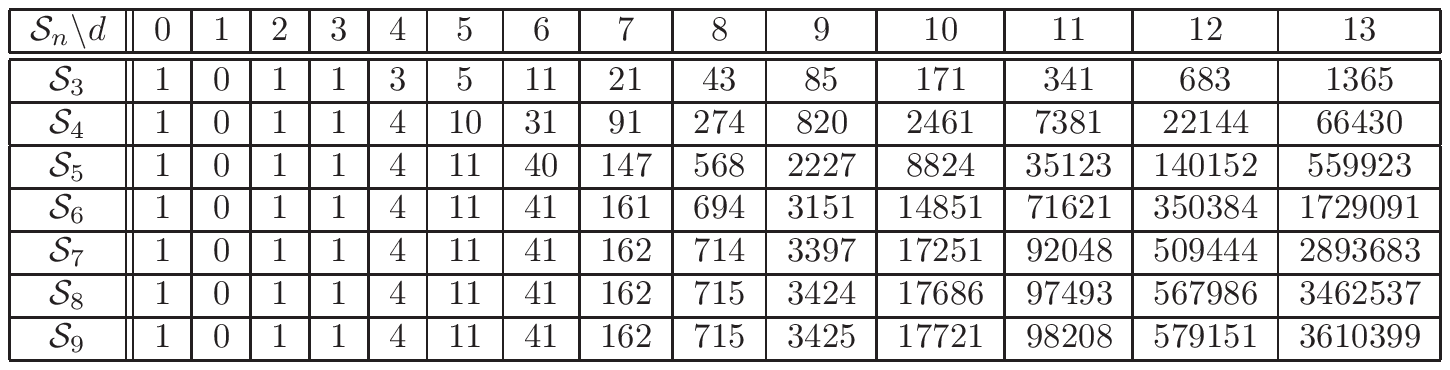}
  \caption{Dimension of $((V^{(n-1,1)})^{\otimes d})^{{\Sn}_n}\simeq\Q\langle Y_{n-1}\rangle_d^{{\Sn}_n}$. Also the number of words of length $d$ which reduce to the identity in $\Gamma({\Sn}_n,\{ (12),(132), (1432), \ldots, (1\,n\,\cdots\,432)\})$.} 
  \label{tab3}
 \end{center}
 \end{table}

\begin{table}[ht!]
\begin{center}
\includegraphics[scale=0.72]{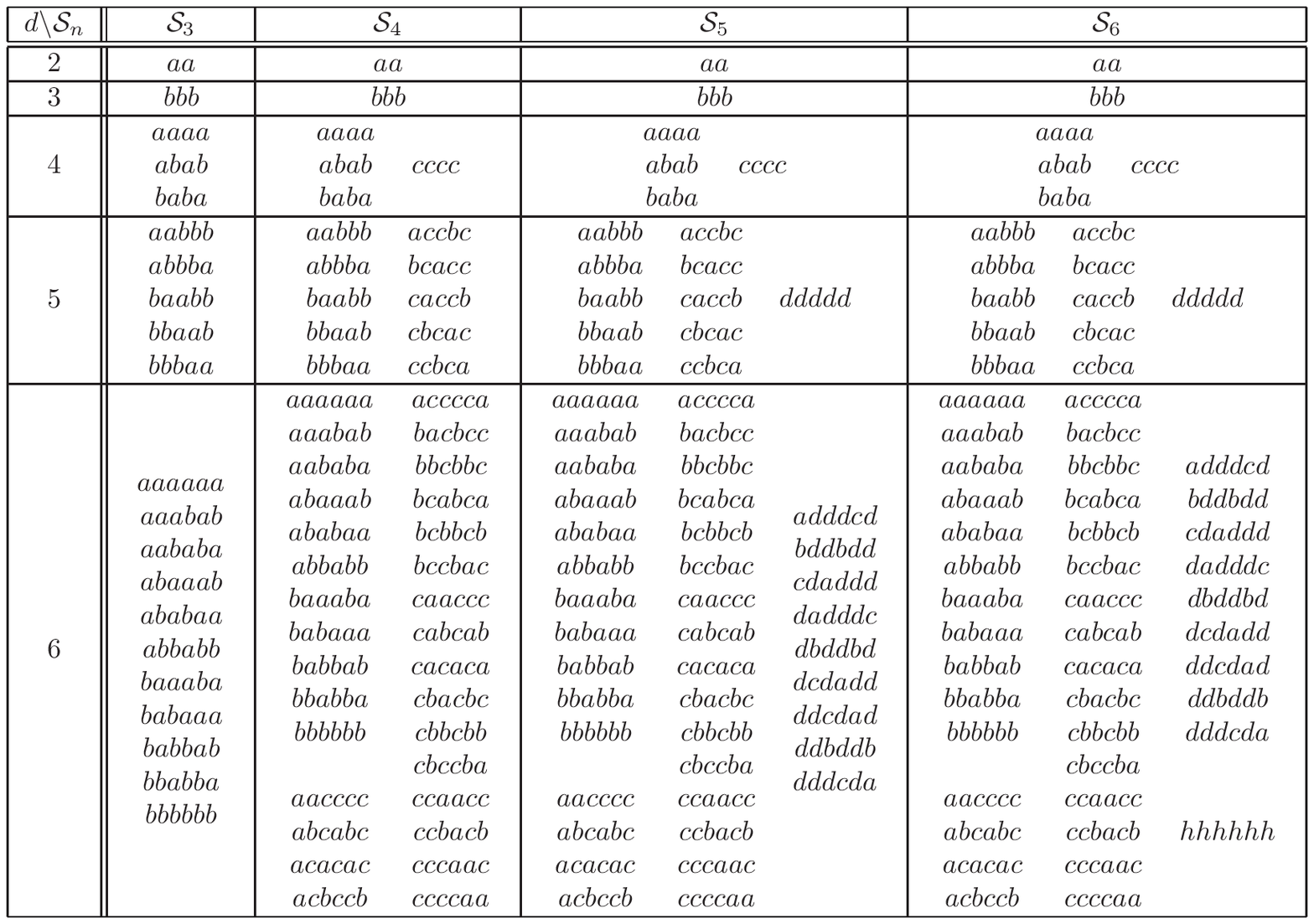}
   \caption{Words of length $d$ in the letters   $a=(12),b=(132), c=(1432), d=(15432), h=(15432)$ which reduce to the identity in $\Gamma({\Sn}_n,\{ (12),(132), (1432), \ldots, (1\,n\,\cdots\,432)\})$.}
     \label{tab4}
 \end{center}
 \end{table}
 
\begin{table}[ht!] 
\label{Aness}
\begin{center}
\includegraphics[scale=0.75]{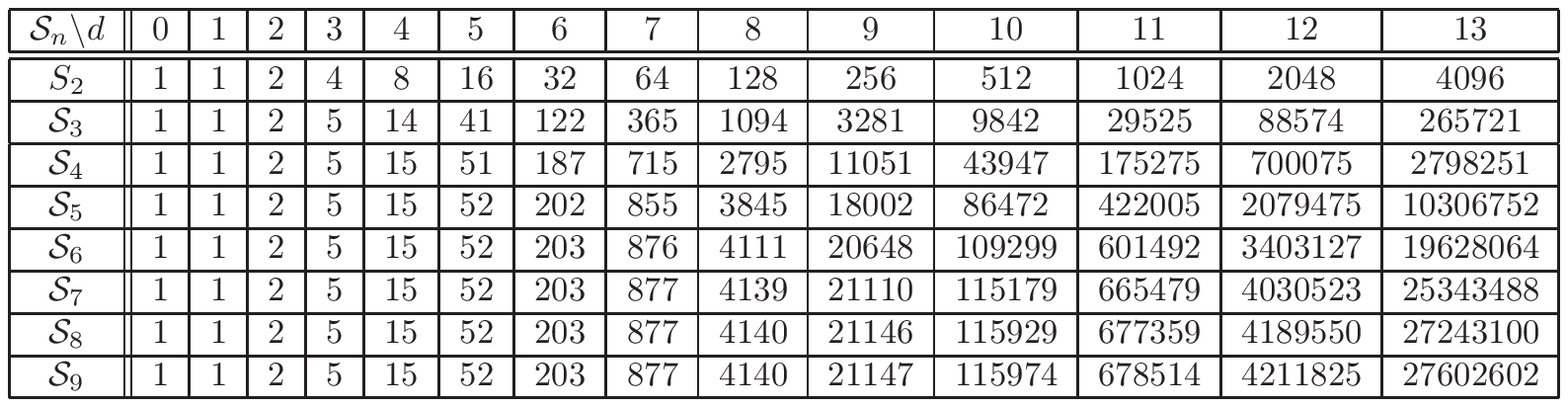}
  \caption{Dimension of  $((V^{(n)}\oplus V^{(n-1,1)})^{\otimes d})^{{\Sn}_n}\simeq\Q\langle X_{n}\rangle_d^{{\Sn}_n}$. Number of words of length $d$ which reduce to the identity in  $\Gamma({\Sn}_n,\{e, (12),(132), (1432), \ldots, (1\,n\,\cdots\,432)\})$. Number of set partitions of $[d]$ into at most $n$ parts.}
    \label{tab1}
 \end{center} 
 \end{table}

\begin{table}[ht!]
\begin{center}
\includegraphics[scale=0.75]{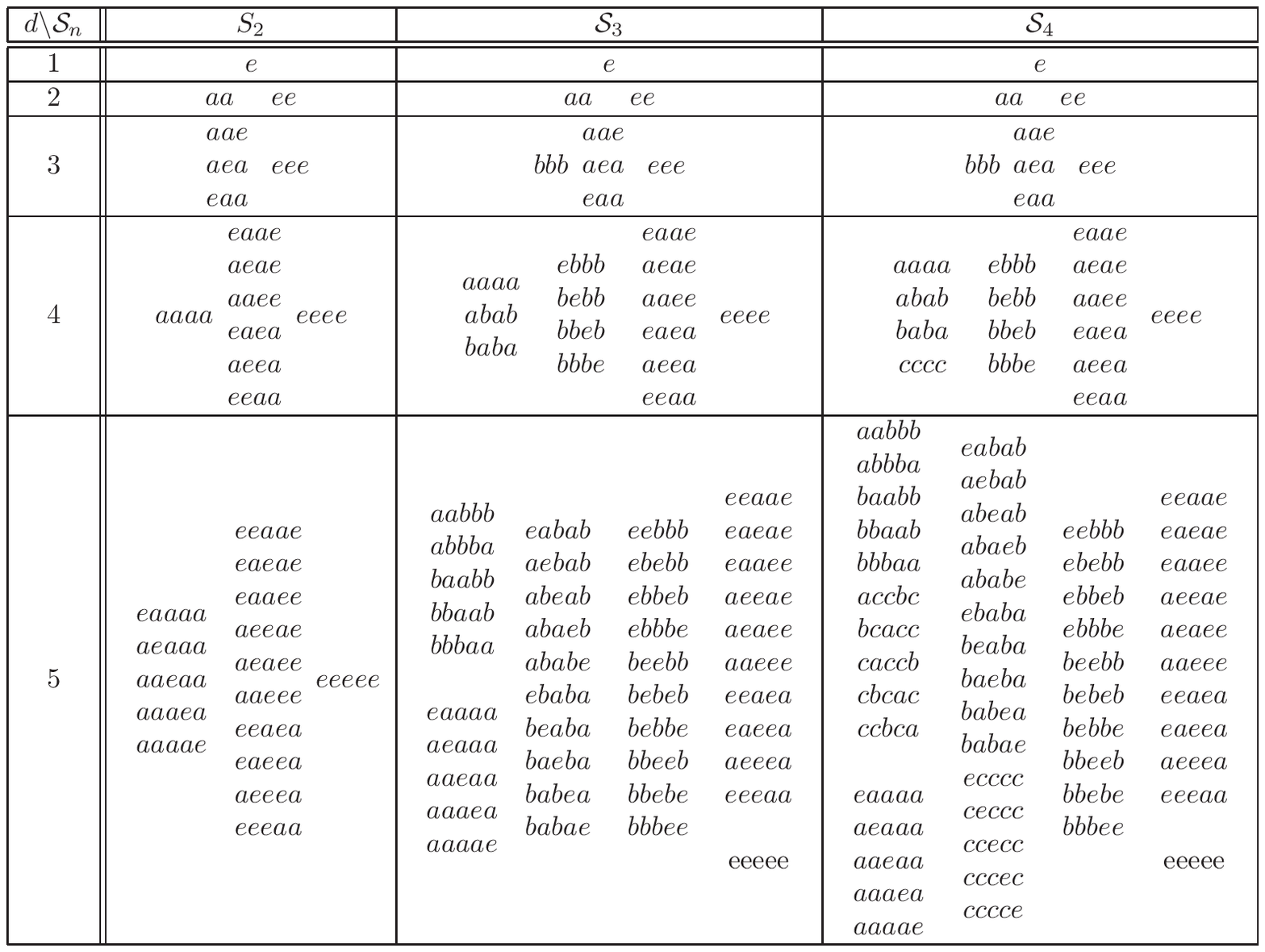}
   \caption{Words of length $d$ in the letters $a=(12)$, $b=(132)$, $c=(1432)$ and $e$ which reduce to the identity in $\Gamma({\Sn}_n,\{e, (12),(132), (1432), \ldots, (1\,n\,\cdots\,432)\})$.}
     \label{tab2}
   \end{center}
 \end{table}

\begin{table}[ht!]
\begin{center}
\includegraphics[scale=0.75]{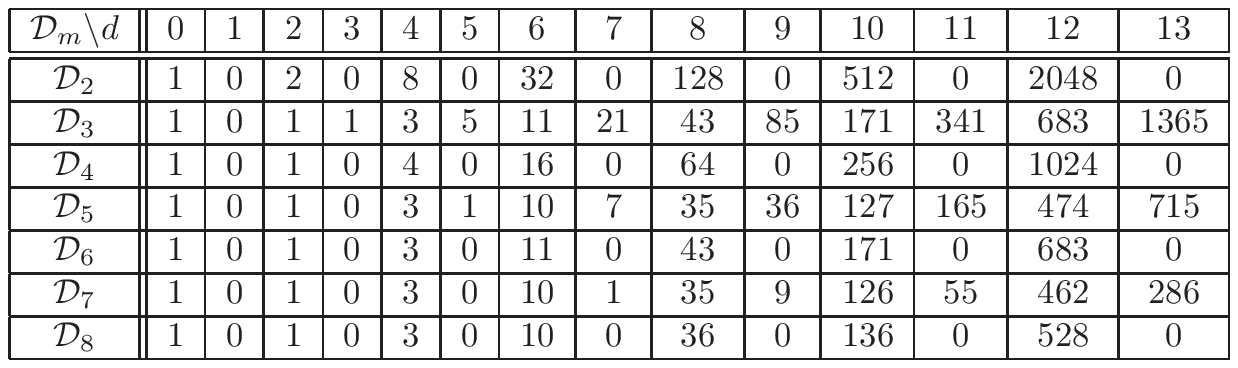}
 \caption{Dimension of  $((V^1)^{\otimes d})^{{\D}_m}\simeq\R\langle x_1,x_2\rangle_d^{{\D}_m}$. Number of words of length $d$ in the letters $r$ and $s$ which reduce to the identity in $\Gamma({\D}_m,\{r,s\})$.}
   \label{tab5}
\end{center}
\end{table}

\begin{table}[ht!]
\begin{center}
\includegraphics[scale=0.75]{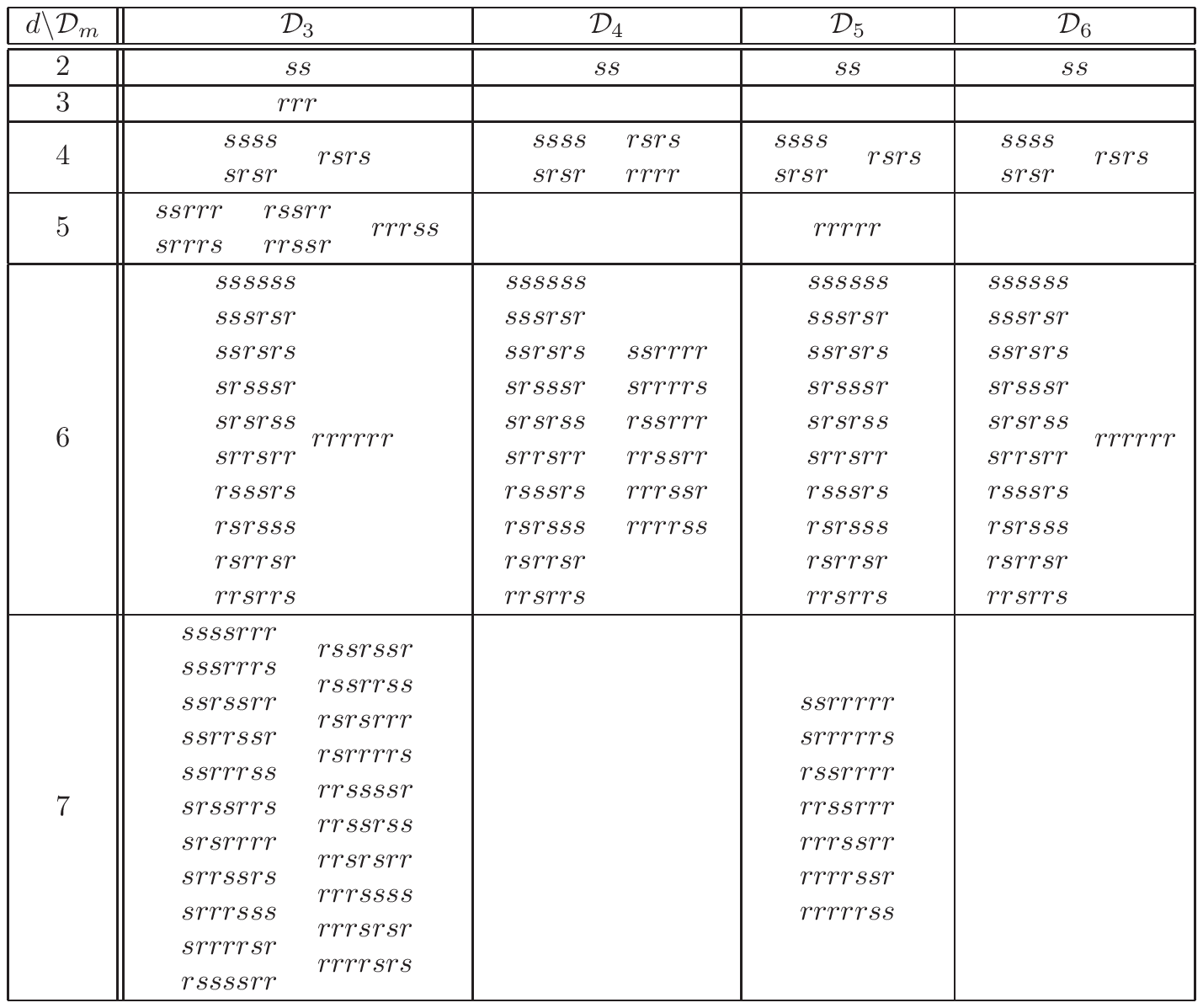}
   \caption{Words of length $d$ in the letters $r$ and $s$  which reduce to the identity in $\Gamma({\D}_m,\{r,s\})$.}
     \label{tab6}
 \end{center}
 \end{table}

\clearpage
\bibliographystyle{plain}
\bibliography{/Users/Anouk/Documents/Articles/BIBLIO}


\subsection*{Acknowledgment}
We would like to thank Andrew Rechnitzer for great help in the proof of Proposition \ref{prop: SerieI2}.
\end{document}